\documentclass[english, a4paper, 11pt,twoside]{amsart}

\usepackage[top=2.5cm, bottom=2.5cm, left=2.5cm, right=2.5cm]{geometry}
	\usepackage{helvet}
	\usepackage[utf8]{inputenc}
	\usepackage[T1]{fontenc}
	\usepackage{enumitem}
	\usepackage{subcaption}
	\usepackage{wrapfig}
	\usepackage{caption}
	\usepackage{multicol}

\usepackage[english]{babel}
	\usepackage{hyphenat}

\usepackage{titlesec}
	\titleformat{name=\section}[hang]{\scshape \large \centering}{\thesection.~}{0em}{}[]
	\titleformat{\subsection}[runin]{\bfseries}{\thesubsection.}{.2em}{}[.]
	\titleformat{\subsubsection}[runin]{\itshape}{\normalfont \thesubsubsection.}{.2em}{}[.]

\usepackage{titletoc}
	\titlecontents{section}[2pc]{\addvspace{1pc}}
	{\contentslabel[\textsc{\thecontentslabel}]{1pc}}
	{}{\dotfill\contentspage}
	[\addvspace{2pt}]
	
	\titlecontents{subsection}[4pc]{}
	{\contentslabel[\subsectionname\thecontentslabel]{2pc}}
	{}{\hfill\contentspage}
	[]

\usepackage{amssymb,amsthm,latexsym}
	\usepackage[fixamsmath]{mathtools} 
	\usepackage{mathabx} 
	\usepackage{MnSymbol}
	\usepackage{tikz-cd}
	 	\usetikzlibrary{decorations.pathmorphing}
\usepackage{graphicx}
\usepackage[dvipsnames]{xcolor}
	\usepackage{colortbl}
	\definecolor{e-mail}{rgb}{0,.40,.80}
	\definecolor{reference}{rgb}{.20,.60,.22}
	\definecolor{mrnumber}{rgb}{.80,.40,0}
	\definecolor{citation}{rgb}{0,.40,.80}
	\definecolor{gris25}{gray}{0.45}
\usepackage{enumitem}
	\setlist[enumerate]{label*=(\roman*)} 

\usepackage[
pdfauthor={Benjamin Collas},
pdftitle={Anabelian perspectives in Galois-Teichmüller theory},
bookmarksnumbered=true,
colorlinks=true,
linkcolor=reference, citecolor=citation, urlcolor=e-mail
]{hyperref}

\usepackage{filecontents}

\usepackage{csquotes}
\usepackage[
bibencoding=utf8,%
giveninits=true,
bibstyle=alphabetic,
citestyle=alphabetic,%
ibidtracker=true,
sorting=nyt,%
language=english,%
isbn=false,%
doi=false,%
eprint=true,
backref=true, 
hyperref=true,
url=false, 
minbibnames=4,
maxbibnames=4,
backend=biber,
]{biblatex}
\AtBeginBibliography{\small\vspace{-.5em}}
\defbibenvironment{bibliography}
{\list
	{\bfseries [\printfield[labelnumberwidth]{labelalpha}\printfield[labelnumberwidth]{extraalpha}]}
	{\setlength{\labelwidth}{0pt} 
		\setlength{\leftmargin}{0pt} 
		\setlength{\labelsep}{4pt}
		\setlength{\itemsep}{4pt}
		\setlength{\parsep}{2pt}} 
	}
{\endlist}
{\item}

\theoremstyle{plain}
	\newtheorem{theo}{Theorem}[section]
	\newtheorem{prop}[theo]{Proposition}

	\newtheorem*{theorem*}{Theorem}

\theoremstyle{definition}
	\newtheorem{defi}[theo]{Definition}

\theoremstyle{remark}
	\newtheorem{rmk}[theo]{Remark}

\numberwithin{equation}{section}

	\newcommand{\MM}{\mathcal{M}}
	
	\newcommand{\QQ}{\mathbb{Q}}
	\newcommand{\ZZ}{\mathbb{Z}}
	\newcommand{\CC}{\mathbb{C}}

	\newcommand{\piet}{\pi_1^{\text{\textup{\'et}}}}

\newcommand{\Aut}{\operatorname{\rm Aut}}
\newcommand{\Isom}{\operatorname{\rm Isom}}
\newcommand{\Out}{\operatorname{\rm Out}}
	\newcommand{\AutNS}[1]{{\rm Aut}_{#1}^\sharp}
	\newcommand{\OutNS}[1]{{\rm Out}_{#1}^\sharp}	
	\newcommand{\GT}{\widehat{GT}}

	\newcommand{\conf}[1]{\operatorname{Conf}_{#1}(X_{0,3})}	
	
	\newcommand{\rO}{\overset{\Out}{\rtimes}}
	
	\newcommand{\BGT}{{\rm BGT}}

\title[Anabelian perspectives in Galois-Teichmüller theory]{Anabelian perspectives in Galois-Teichmüller theory \\[0.5em]
	{\normalfont \footnotesize In honor of Professor Hiroaki \textsc{Nakamura}'s $60$th birthday}}
	\date{Version of \today}

\author[B.~Collas]{Benjamin~Collas}
	\address{Research Institute for Mathematical Sciences, Kyoto University, Kyoto 606-8502, Japan}
	\email{bcollas@kurims.kyoto-u.ac.jp}

\subjclass[2020]{Primary: 14H30, 11G99; Secondary: 14H10, 20F36, 20F34, 11R32, 14F35}
\keywords{Arithmetic homotopy geometry, combinatorial and mono-anabelian geometry, Grothendieck-Teichmüller theory, moduli spaces of curves, braid groups and configuration spaces of curves.}
	\thanks{The author would like to thank Shinichi \textsc{Mochizuki}, Yuichiro \textsc{Hoshi}, and Shota \textsc{Tsujimura} for numerous exchanges on combinatorial anabelian geometry, and Arata \textsc{Minamide} for enlightening talks on ``GT anabelian''. The author thanks the anonymous referee for a careful reading and numerous corrections that improved the precision of this survey. This manuscript is part of the France-Japan AHGT international research network supported by the Research Institute for Mathematical Sciences of Kyoto University and CNRS. The author acknowledges the support of the International Center for Next-Generation Geometry, a center affiliated with the Research Institute for Mathematical Sciences located in Kyoto University.}

\begin{document}

\newgeometry{bottom=1.5cm, top=2cm, left=2.5cm, right=2.5cm}
	\begin{abstract}
		By exploiting the arithmetic homotopy of the moduli spaces of curves, Galois-Teichmüller theory stands at the interface of braid-mapping class groups and of anabelian geometry. Starting from the classical braid-theoretic construction of the Grothendieck–Teichmüller group $\GT$, we review  how anabelian geometry -- beginning with the foundational work of Nakamura -- provides the arithmetic mechanisms underlying its definition. We then explain how the combinatorial anabelian geometry developed by Hoshi and Mochizuki recasts these constructions within a purely group-theoretic and algorithmic framework. In particular, we describe how the group $\GT$ emerges as an anabelian object and how, once freed from an auxiliary or artificially imposed container, the anabelian algorithms yield a combinatorial reconstruction of the absolute Galois group of rational numbers. The perspective developed here highlights a conceptual shift from explicit braid-theoretic computations to functorial and algorithmic forms of anabelian reconstruction.
	\end{abstract}
	
	\maketitle
	
	\vspace{-1cm}
	
	\tableofcontents
	\clearpage
	
\restoregeometry	
	
	\section*{Introduction}
		
		Galois-Teichmüller theory exploits, via the étale fundamental groups of moduli spaces of curves, their arithmetic homotopy structure in order to obtain a combinatorial description -- originally formulated in terms of explicit representations of braid groups and mapping class groups -- of the absolute Galois group of the rational numbers. The conceptual foundation of this approach lies in the anabelian nature of the moduli spaces of curves. This anabelian perspective provides a universal anchor for the homotopy-theoretic approach, in the sense that homotopy faithfully encodes both the arithmetic and the geometry of these spaces.
		
		\medskip
		
		From this point of view, the aim is to construct intrinsic group-theoretic objects -- of which the Grothendieck-Teichmüller group may be regarded as an intermediate test object -- or, in more recent developments such as combinatorial anabelian geometry, to formulate functorial group-theoretic algorithms that capture and reconstruct these arithmetic homotopy properties.
		
		\medskip
		
		In the first section, we review the classical Grothendieck-Teichmüller program, recalling the arithmetic origin of its braid-theoretic combinatoric, and introduce the anabelian geometry principles that will guide the remainder of this manuscript. The second section explains how the anabelian properties of $\MM_{0,m}$ lead naturally to the notion of quasi-special automorphisms, which once organized within the group $\GT$, then emerges as the stable outer automorphism group of the genus-zero modular tower. The third section develops the techniques of combinatorial anabelian geometry that allow one to remove the auxiliary cuspidal arithmetic conditions and to obtain a purely group-theoretic anabelian characterization of $\GT$. The final section completes the transition from computational braid-theoretic methods to algorithmic anabelian reconstruction, culminating in a combinatorial description of the absolute Galois group $G_\QQ$. 
    
    \medskip

    Let us emphasize the movement that takes place: in \S~\ref{sub:GTtower} the group $\GT$ is introduced as a \emph{container} that serves to organize the braid-theoretic computations; in \S~\ref{sec:GTAnab} it becomes an anabelian object in its own, so that in \S~\ref{sec:BGT} it is finally \emph{dispensed with} for the reconstruction of $G_\QQ$ to rely on purely anabelian properties of Galois-Teichmüller theory, see Rem.~\ref{rm:freeGT}.

		\bigskip
		
		{\small 
		※ This manuscript is based on a series of lectures delivered in recent years, including the ``Inter-universal Teichmüller Theory Summit 2025'' (Kyoto, JP), the ``Journées galoisiennes d'hiver'' (Caen FR, 2024), the ``22nd Affine Algebraic Geometry Meeting'' (Niigata JP, 2024), and ``the Grothendieck's approach to Mathematics'' conference (Chapman University US, 2022). Precise references are provided to enable readers to turn their experience into an active mathematical practice.
		}
		
	\section{Arithmetic homotopy foundation of Galois-Teichmüller theory}

		\subsection{The geometric Galois action on $\MM_{0,[4]}$}\label{sub:GGA}
			For $X$ a quasi-projective and geometrically connected variety over $\QQ$, and $*\colon {\rm Spec }\ \bar{\QQ}\to X$ a geometric point, \cite{SGA103} provides a homotopy fundamental exact sequence (FES), then -- even equivalently when $Z(\piet(X\times \bar{\QQ},*))=\{1\}$ -- an outer Galois action as below (RHS):
			\[\begin{tikzcd}[column sep=1em]
				1 \ar[r] & \piet(X\times \bar{\QQ},*) \ar[r]&  \piet(X,*) \ar[r] & G_\QQ \ar[r]&1
			\end{tikzcd}\quad \text{defines} \quad G_\QQ \to \operatorname{Out}(\piet(X\times \bar{\QQ},*)).
			\]
			Once an embedding $\bar{\QQ}\hookrightarrow \CC$ is fixed, one further identifies $\widehat{\pi}_1^{\rm top}(X(\CC)^{\rm an},*)\simeq\piet(X\times \bar{\QQ},*)$ to bring specific topological combinatorics in the study of the absolute Galois group $G_\QQ$.
			 
			\medskip
			
			Galois-Teichmüller theory, following Grothendieck's insight in \cite{GroEsq97}, deals with considering the moduli spaces of curves of genus $g$ with $m$-marked points $X=\MM_{g,[m]}$ and exploits
			\begin{enumerate}
				\item a $G_\QQ$-invariant stratification of $\MM_{g,[m]}$ by lower dimensional $(g',m')$-strata;
				\item an étale-analytic topology comparison;
				\item a mapping class group or braid group combinatorial description of $\widehat{\pi}_1^{\rm top}(X(\CC))$
			\end{enumerate}
			to obtain \emph{a combinatorial description of $G_\QQ$.}
			
			\subsubsection{``La droite projective moins $3$ points"}\label{sub:P3pts}
\begin{wrapfigure}[7]{r}{5cm}
					\vspace*{-1.2em}
					\centering
					$\begin{tikzcd}[column sep=-2em, row sep=.7em]
						&\bar{\QQ}\{\{t\}\}\simeq \bar{\QQ}\{\{1-t\}\} & \\
						M_{\vec{01}}\ar[rr,"p","f_\sigma" ',squiggly] \ar[ru]& &M_{\vec{10}} \ar[lu]\\
						& \bar{\QQ}(t)\ar[ru, hook] \ar[lu,hook']& 
					\end{tikzcd}$
					\captionof{figure}{Analytic transport of the étale $\chi(\sigma)$}\label{fig:Qp}
				\end{wrapfigure}
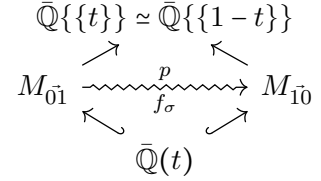

				We briefly recall, following \cite{Iha94GT}, the case of $\MM_{0,[4]}=\MM_{0,4}/S_4$, where $\MM_{0,4}(\CC)\simeq \mathbb{P}^1\setminus\{0,1,\infty\}$ and $\pi_1^{\rm top}(\MM_{0,4}(\CC),*)\simeq \mathbb{F}_2$, the free group of rank $2$, has two homotopic generators $x$ and $y$, respectively a loop around $0$ and $1$. Given $\sigma\in G_\QQ$, the outer Galois action is computed as:\\
				\begin{minipage}{.5\textwidth}
					\begin{align*}
						\sigma.x&=x^{\chi(\sigma)}\\
						\sigma.y&=f_{\sigma}^{-1}.y^{\chi(\sigma)}.f_{\sigma}
					\end{align*}
				where $\chi(\sigma)\in\widehat{\ZZ}^\times$ is the cyclotomic character resulting from the monodromy around $0$, and $f_{\sigma}\in \widehat{\mathbb{F}}_2'$ appears by transporting the \emph{étale} Galois action via an \emph{analytic} path $p$ from $0$ to $1$  -- see Fig.~~\ref{fig:xyp}-\ref{fig:Qp}.	
				\end{minipage}\hfill
				\begin{minipage}{.45\textwidth}\centering
					\includegraphics[width=\linewidth]{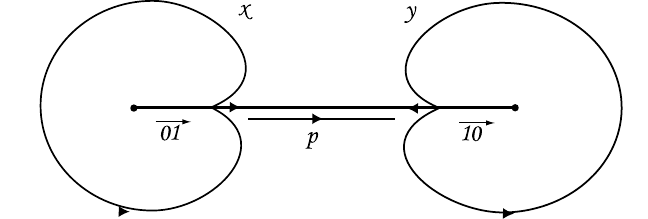}
					\captionof{figure}{Loops, path, and tangential base points in $\MM_{0,4}$}\label{fig:xyp}
				\end{minipage}
								
				\medskip
				
				Note that, in the above, the geometric base point $*$ is indeed replaced by \emph{Deligne-Ihara's tangential base point $\vec{01}$} -- formally a $\QQ$-formal neighbourhood of $0$ and a fixed parameter $t$ of $\QQ[[t]]$ -- which also provides a section to the (FES), and thus a lift
				
				\medskip
				 
				{\centering
					$\varphi_{\vec{01}}\colon G_\QQ\to \Aut(\widehat{\pi}_1^{\rm top}(\MM_{0,4}(\CC)^{\rm an},\vec{01}))$
				\par}
				
				\medskip	
				
				to the canonical outer Galois action. Denoting $M_{\vec{01}}$ (resp. $M_{\vec{10}}$) the maximum Galois extensions unramified outside $\{0,1,\infty\}$ of the function field of $\MM_{0,4}$ in the Puiseux algebraic closure corresponding to $\vec{01}=\bar{\QQ}[[t]]$ (resp. to $\vec{10}=\bar{\QQ}[[1-t]]$), Fig.~\ref{fig:Qp} illustrates the analytic transport of the cyclotomic $G_\QQ$-action from $\vec{01}$ (at $0$) to $\vec{10}$ (at $1$). 
				
				\medskip

				The $G_\QQ$-action, once $\widehat{\ZZ}^\times\times \widehat{\mathbb{F}}_2'$ is endowed with the composition law corresponding to the $G_\QQ$-composition given by $\varphi_{\vec{01}}$, provides in turn a \emph{homomorphism} $G_\QQ\to \widehat{\ZZ}^\times\times \widehat{\mathbb{F}}_2'$.
				
			\subsubsection{Homotopy symmetries and $\GT$}\label{sub:Hsym} 
				To get rid of any non-canonical choice in the previous construction, we keep track of the $S_4$-symmetries of $\MM_{0,4}$, by further \textbf{(a)}~endowing $\widehat{\pi}_1^{\rm top}(\MM_{0,4}(\CC),\vec{01})$ with a structure of groupoid by taking the set $\mathcal{B}$ of all the six tangential base points of $\MM_{0,4}$, and \textbf{(b)}~implementing the $S_4$-symmetries in terms of certain homotopy relations.
				\begin{itemize}
					\item[(I)] $f_{\sigma}(x,y).f_{\sigma}(y,x)=1$
					\item[(II)] $f_{\sigma}(z,x)z^m.f_{\sigma}(y,z)y^m.f_{\sigma}(x,y)x^m=1$ 
				\end{itemize}
				where $m=(\chi(\sigma)-1)/2$ and $xyz=1$, and with a third one denoted (III) that is of similar nature with respect to $\MM_{0,5}$ -- see \cite{Iha94GT}.
				
				\medskip
				
				This results in a purely combinatorial description of $G_\QQ$ -- see Ihara \cite{Iha83} \S~3-4 and \cite{Iha94GT}, and also \cite{Dri91} for a similar approach in terms of quasi-triangular quasi-Hopf algebra and monoidal braided categories.
				\begin{theo}
					There is a Grothendieck-Teichmüller group $\GT$ of group elements $(\lambda,f)\in \widehat{\ZZ}^\times\times \widehat{\mathbb{F}}_2'$ satisfying relations (I), (II), and (III) that acts as $G_\QQ \hookrightarrow \Aut(\widehat{\pi}_1^{\rm top}(\MM_{0,4}(\CC),\vec{01}))$ via $\varphi_{\vec{01}}$, and contains $G_\QQ\hookrightarrow \GT$.
				\end{theo}
				
				\begin{wrapfigure}[7]{r}{5cm} 
        \vspace*{-1.5em}
        \centering
        $\begin{tikzcd}
          G_\QQ \ar[r,hook] \ar[rd,hook]& \Out(\widehat{\Gamma}_{0,[m]})\\
          & \GT \ar[u, hook]
        \end{tikzcd}$
        \captionof{figure}{Arithmetic, geometric, and group-theoretic triangle}\label{fig:AGtr}
      \end{wrapfigure}
				The homotopies relations are obtained by considering $S_4\simeq V_4\rtimes S_3$ acting on $\MM_{0,4}$ via $S_3=\langle \theta,\omega\rangle$ and the order two (resp. order three) $G_\QQ$-invariant automorphisms $\theta\colon t\to 1-t$  (resp. $\omega\colon t\to 1/(1-t)$), which respectively provide \textbf{(I})~$\theta(p).p=1$ and \textbf{(II)}~$\omega^2(q).\omega(q).q=1$ for well chosen paths $p$ and $q$ -- see \cite{Iha83} for exact definition and \cite{EL97} some explicit and enlightening computations. The Galois injectivity famously follows from Belyi's theorem.

		 \subsection{Higher moduli spaces and braid arithmetic}\label{sub:higher} 
		 	Let us fix $g=0$ and consider the higher dimensional spaces $\MM_{0,[m]}$ for $m\geq 4$. The context and results of \S~\ref{sub:GGA} extend via the combinatoric of braid groups and the sophisticated notion of tangential structures.
		 	
		 	\subsubsection{Tangential structures and braid groups}\label{sub:HdimAH} 
		 		The previous context of \S~\ref{sub:P3pts} and \S~\ref{sub:Hsym} extends as follows:
		 		\begin{enumerate}
		 			\item The orbifold fundamental group $\pi_1^{\rm orb}(\MM_{0,[m]}\times \bar{\QQ})\simeq \Gamma_{0,[m]}$ is the full mapping class group of diffeomorphisms, that is identified with \emph{Artin braid groups} $\Gamma_{0,[m]}\simeq B_m/\langle H_m,Z_m\rangle$, quotiented by the Hurwitz and the  center relations;
		 			\item A \emph{tangential structure} on $\MM_{0,[m]}$ is given by a system of parameters $\mathbf{t}=\{t_1,\dots t_{m-3}\}$ of a $\QQ$-formal neighbourhood\footnote{Geometrically, each $t_i$ corresponds to a simple closed loop on a pant decomposition of the Riemann surface of type $(g,m)$ that degenerates to a singular point on $X_0$.} in $\MM_{0,[m]}$ of a genus $0$ stable curve $X_0$ over $\QQ$ in the Deligne-Mumford compactification $\overline{\MM}^{DM}_{0,[m]}$.
		 			\item The $\GT$-action on braid groups $\widehat{B}_m$ is given by Drinfel'd via quasi-triangular quasi-Hopf algebra \cite{Dri91} -- see Eq.~\eqref{eq:FormGT} below.
		 			\item \emph{The étale fundamental group} associated to a tangential structure is given by Grothendieck-Murre's $\pi_1^{\mathcal D}(\MM_{g,[m]}, \mathbf{t})$ of \cite{GM71}, that classifies étale covers unramified outside a normal crossing divisor $\mathcal{D}$ supporting the curve $X_0$.
		 		\end{enumerate}
		 		
		 		\medskip
		 		
		 		We further refer to \cite{IN97max} for the definition of tangential structures via $\mathbb{P}^1_{0,1,\infty}$-diagrams (see also Fig.~\ref{fig:NakTg}), and  to \cite{CM23} \S~3.1 for a more intrinsic stack definition (see also Rem.~3.6 ibid.).
		 	 
		 	 	\medskip
		 	 	
		 	 	{\centering
		 	 	\includegraphics[width=.7\linewidth]{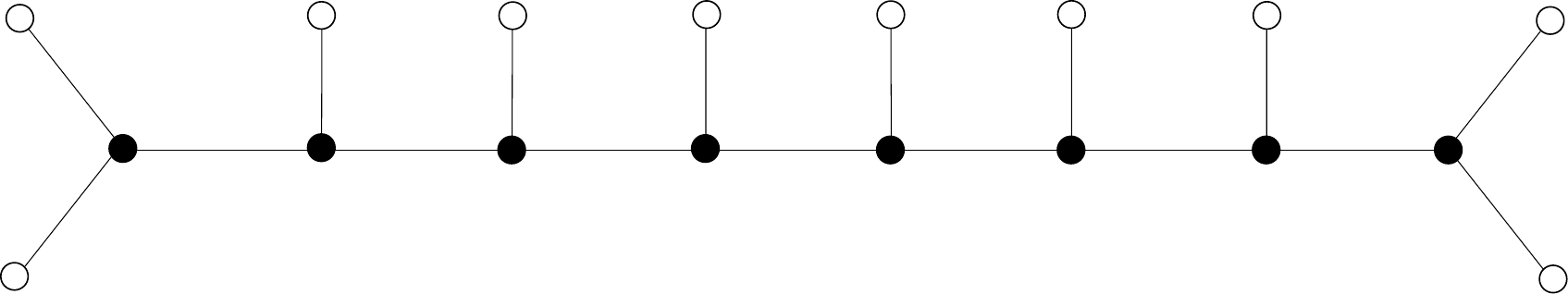}
		 	 	\captionof{figure}{Dual graph of the stable curve $X_0$ serving as tangential structure.}\label{fig:NakTg}
		 	 	\par}
		 	 	
		 		\medskip
		 		
		 		Taking Nakamura's tangential structure of \cite{Nak97tg} Eq.~(4.2) associated to Fig.~\ref{fig:NakTg} or, equivalently, following the braid-computations and blow-up techniques of configuration spaces of \cite{IM95}, one further obtains. 
		 		\begin{theo}\label{th:GTBn}
		 			For $m\geq 4$, there exists a tangential $G_\QQ$-action on $\piet(\MM_{0,[m]}\times\bar{\QQ})$ given on the Artin generators $\sigma_i$, $i=1,\dots, m-1$, of $\widehat{B}_m$ by
		 			\begin{equation}\label{eq:FormGT}
				 		\sigma.\sigma_i=f_{\sigma_i}(y_i,\sigma_i^2)^{-1}.\sigma_i^{\chi(\sigma)}.f_{\sigma_i}(y_i,\sigma_i^2),\quad \text{ where }y_i=\sigma_i\sigma_{i-1}\dots \sigma_1.\sigma_1\dots 	\sigma_{i-1}\sigma_i,
		 			\end{equation}
		 			with $y_1\coloneq 1$, that is compatible with the $\GT$-action on $\widehat{B}_m$.
		 		\end{theo}
		 		
		 		It results from the explicit definition of the $\GT$-action on $\widehat{B}_m$ above, that the $\GT$-action stabilizes both the Hurwitz relation $H_m$ and the center relation $Z_m$ so that it passes to the quotient $\widehat{\Gamma}_{0,[m]}\simeq \widehat{B}_m/\langle H_m,Z_m\rangle$. The group $\GT$ also acts on the pure braid group $\widehat{P}_m<\widehat{B}_m$ as a consequence of Dyer-Grossman's characteristic property for $P_m$. 
		 		
		 		\medskip
		 		
		 		It may be of interest to note that, while an expression of the form $\sigma.\sigma_i=A_{\sigma,\sigma_i}^{-1}. \sigma_i^{\chi(\sigma)}.A_{\sigma,\sigma_i}$ can readily be obtained from the action of \S~\ref{sub:P3pts} -- since every generator $\sigma_i$ is conjugated to $\sigma_1$ -- the higher dimensional action above is entirely defined \emph{in terms of the dimension one and two moduli spaces of genus zero}, that is of $\MM_{0,[4]}$ and $\MM_{0,[5]}$ -- which follows from Grothendieck's $2$-level guiding principle, see \cite{GroEsq97}.

				\subsubsection{In higher genus...}
					Let us briefly mention that, in higher genus $g\geq 0$, there also exists a Grothendieck\hyp{}Teichmüller group $\rm I\!\Gamma$: it originates from a certain Thurston-Hatcher two-dimensional complex of pants decompositions on Riemann surfaces of type $(g,m)$, and relies on explicit presentations of the mapping class group $\Gamma_{g,m}$, so that one has an action
					\[
						{\rm I\!\Gamma}\longrightarrow\Out(\widehat{\Gamma}_{g,[m]}) \text{ in a $G_\QQ$-compatible way, with }G_\QQ\leq {\rm I\!\Gamma} \leq \GT\quad 
					\]
					We refer to \cite{Nakamura00} for the original definition and to \cite{Colg1} for its complete action on the full mapping class group $\widehat{\Gamma}_{g,[m]}$.

		\subsection{From classical anabelian geometry to combinatorial anabelian geometry} 
			The anabelian mechanisms involved in this manuscript rely on a contemporary shift from explicit computational methods to algorithmic reconstruction, articulated through mono-anabelian and combinatorial anabelian frameworks.
			
			\subsubsection{Anabelian arithmetic geometry principles}\label{subsub:AGprinc}
				Anabelian arithmetic geometry deals with the reverse ``functorial'' process which, to a space $X$ over a field $K$, attaches a homotopy fundamental exact sequence (FES) as in \S~\ref{sub:GGA}. 
				
				\medskip
				
				Considering certain classes of varieties  $\mathfrak{X}_K$, anabelian arithmetic geometry is motivated by Grothendieck's 1983 foundational insight\footnote{On contrary to Grothendieck's original statement, the class of Siegel modular varieties is not anabelian, see~\cite{IN97}.} \cite{GrFal97} -- where the subscript $G_K$ denotes ``$G_K$-equivariant morphisms'':
				
				\medskip
				
				{\centering
				\begin{minipage}{.9\textwidth}\itshape
					\textbf{(GC)~The Fundamental Conjecture of Anabelian Geometry\footnote{The emphasis on arithmetic aspects of algebraic geometry is subsequent to Grothendieck and follows Deligne's IHES ``séminaire de géométrie arithmétique'' that follows on dominant cohomological approach; hence our use, in the present manuscript, of the term ``anabelian arithmetic geometry'' or even ``arithmetic homotopy geometry''.}} The isomorphism class of a variety $X$ in $\mathfrak{X}_K$ can be reconstructed from its étale fundamental group, that is, for $X$ and $Y$ in $\mathfrak{X}_K$, the application
					\[
					\Isom_K(X,Y)\to \Isom_{G_K}(\piet(X),\piet(Y))/\operatorname{Inn}\piet(Y\times\bar{K})
					\]
					is a bijection. A variety of a class $\mathfrak{X}_K$ that satisfies (GC) is called ``anabelian''.
				\end{minipage}
				\par }
				
				\medskip
				
				In dimension zero, \emph{a number field} $X=\operatorname{Spec} K$ is anabelian by the 1970s Neukirch-Uchida theorem, see for example \cite{NSW08} Theorem~12.2.1. In dimension one, it further results from the incremental efforts of Nakamura \cite{Nak90} (genus $0$, over number fields) and Tamagawa \cite{Tam97} (affine curves, over number fields), then Mochizuki \cite{MOC99} (proper curves) that the class of hyperbolic curves over a sub-$p$-adic field\footnote{A class which contains the number fields as well as the \emph{mixed-characteristic local fields} (MLF), that is the finite extensions of $\QQ_p$.} $K$ is anabelian. In higher dimension, the moduli spaces $\MM_{0,m}$ are anabelian -- see \S~\ref{subsub:AnabNak}.
				
				\medskip
				
				Noticing, however, that (GC) fails for $X=\operatorname{Spec} K$ an MLF (see \cite{NSW08} above Th.~12.2.7), that is ``there exist isomorphisms of Galois groups of MLF that don't originate from scheme theory'', resulted in a shift of anabelian arithmetic geometry toward a more \emph{algorithmic and absolute approach} -- absolute in the sense that a scheme $X$ is not endowed with its structural morphism $X\to\operatorname{Spec} K$, and mono-anabelian in the sense that it involves a unique object -- see \cite{AbsTopI} \S~Introduction.
				
				\medskip
				
				{\centering
					\begin{minipage}{.9\textwidth}\itshape
						\textbf{(AlgAG)~Algorithmic Anabelian Geometry.} Establish a functorial topological group-theoretical algorithm, which, from $\piet(X)$ generates, beyond the class of $X$ itself, discrete data from geometric or arithmetic origin.
					\end{minipage}
					\par }
				
				\medskip

				The ``functoriality'' requirement excludes non-canonical choices (and thus may introduce the consideration of orbits of objects); ``topological'' refers to compatibility with the profinite topology (note that the $p$-adic topology of $G_{\QQ_p}$ is not group-theoretic compatible). 
				
				In the case of the Galois group of MLF, we refer to \cite{Hos21} for examples of algorithms and of data which rely on classical constructions of local class field theory -- see also perspectives in \cite{CMY25} \S~1.2 and \S~2.2.
				
				\medskip
				
				\pagebreak

				Over the last $20$ years, this \emph{algorithmic anabelian geometry}:
				\begin{enumerate}
					\item has matured in \textbf{a series of anabelian algorithms}, for MLF \cite{Hos21}, and for hyperbolic curves of strict Belyi type (see Fig.~\ref{Fig:StcB}) to reconstruct their function field (resp. their base field) -- see  \cite{AbsTopIII} Cor.~1.11~(iii) (resp. Th.~1.9 (d)-(e)) -- as well as various discrete data such as genus, number of marked points, and cuspidal inertia groups;
					\item has put \textbf{new structures into light} which are comparatively \emph{weaker} to the ring or field ones -- the structure of \emph{multiplicative monoid}, which is one of the usual outputs of (AlgAG), and
					\item developed \textbf{a mono-analytic (MA) vs arithmetic holomorphic (AH) philosophy}. From the point of view of anabelian algorithms, (MA) objects are one-dimensional and strictly related to the multiplicative monoidal structure, while (AH) objects are anabelian global objects; the group-theoretic automorphisms of the (AH) objects play the role of \emph{a variation of analytic structure }for the (MA) objects -- in a similar way that $f_{\sigma}$ appears in \S~\ref{sub:GGA} -- see \cite{AbsTopIII} \S~I3 .
				\end{enumerate} 
				
				\medskip
				
				In the sense of (AlgAG) and the (MA)-(AH) guiding philosophy, the combinatorial computations of \S~\ref{sub:HdimAH} in terms of braid groups and tangential structure take the form of a \emph{Combinatorial anabelian geometry} (CmbAG) -- see \S~\ref{sub:whatisCbAn} and  \S~\ref{sec:GTAnab}. Despite \emph{a radical shift of perspective}, (CmbAG) is not disconnected from the original approach of \S~\ref{sub:GGA} and (GC); for example, the rigid data of the image of $G_\QQ$ in $\Out(\piet(\MM_{0,m}\times\bar{\QQ}))$ and the Galois centralizer $Z_{\Out(\piet(\MM_{0,m}\times\bar{\QQ}))}(G_\QQ)$, or its group-theoretic variant $S_m$ (see below Eq.~\eqref{eq:GTanab}), are of equal important interest to both.
				
			\subsubsection{What is combinatorial anabelian geometry?}\label{sub:whatisCbAn}
				In a few words
			
				\medskip
				
				{\centering
					\begin{minipage}{.9\linewidth}\itshape
						\textbf{Combinatorial anabelian geometry (CombAG)} can be seen as a purely group-theoretic version, in the sense of (AlgAG), of the combinatorial $G_\QQ$-action properties on $\piet(\MM_{g,m}\times\bar{\QQ})$ previously obtained in terms of braid groups/Lie algebra/mapping class groups and $\mathbb{P}_{0,1,\infty}^1$-diagrams/tangential structures computations,
					\end{minipage}
					\par}
				
				\medskip
				
				\noindent 
				that we encountered in \S~\ref{sub:GGA}-\S~\ref{sub:higher}. As such, it benefits from data of a more canonical and universal nature (for the non-homotopic arithmetic geometer, it can also be seen as a sort of homotopy Picard-Lefschetz theory with finer arithmetic and geometric insights.)
				
				\medskip
			
				More precisely, consider a formal neighborhood of a stable curve $X_0\in \bar{\MM}_{g,m}(\CC)$, that is a hyperbolic $(g,m)$-curve $X$ over $\mathcal{O}_k=\CC[[t]]$ which is generically smooth and has (singular) stable fiber $X_0$ over the residue field of $\mathcal{O}_k$. The usual homotopy (FES) then becomes, see (LHS) below
				\[\begin{tikzcd}[column sep=1em]
					& \widehat{\Gamma}_{g,m} \arrow[d, phantom, "\simeq",sloped]&  & \widehat{\ZZ} & \\[-1em]
					1 \ar[r]& \piet(X\times\bar{k}) \ar[r]& \piet(X)   \ar[r]& G_k\ar[r] \arrow[u, phantom, "\simeq",sloped]& 1
				\end{tikzcd}
				\begin{tikzcd}[column sep=1em]
					 & \\[1em]
					\text{defines}\quad\widehat{\mathbb{\ZZ}} \ar[r]&  \Out(\widehat{\Gamma}_{g,m})
				\end{tikzcd}
				\]
				where \emph{the resulting outer action} of the (RHS) above captures \emph{the monodromy variation on the smooth fiber} as given under the action of a Dehn twist on a Riemann surface of type $(g,m)$. In terms of anabelian reconstruction, (CombAG) deals with the characterization, at the level of étale fundamental groups, of the dual graph (vertices and edges) of the special fiber attached to $X_0$ -- referred to as a \emph{combinatorial Grothendieck conjecture-type}.
				
				\medskip
				
				Considering $X_0$ a Pointed Stable hyperbolic Curve of type $(g,m)$ as above (PSC), one can attach, via its dual semi-graph $\mathbb{G}_{X_0}$, \emph{a semi-graph of anabelioids (of PSC-type) that is a Galois category $\mathcal{G}_{X_0}$}, obtained by gluing at each connected component various Galois categories -- see \cite{CombGC}. 
				
				\pagebreak
				
				More precisely:
				
				\smallskip
				 
				\begin{minipage}{.5\textwidth}\centering
					\includegraphics[width=.8\linewidth]{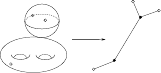}
					\captionof{figure}{From stable curve to dual graph}
				\end{minipage}
				\hfill
				\begin{minipage}{.45\textwidth}
					\begin{itemize}[leftmargin=1em]
						\item At a vertex $v\in \mathbb{G}_{X_0}$ (resp. an edge), $\mathcal{G}_v$ is the category of finite étale covers over $X_0^v$ (resp. over the pointed disk);
						
						\item Each edge $v\rightsquigarrow w$ gives rise to an exact functor $\mathcal{G}_v\to\mathcal{G}_w$.
					\end{itemize}
					One can further form \emph{the PSC-fundamental group $\Pi_\mathcal{G}$} of $\mathcal{G}$, which, by definition contains the verticial subgroups $\Pi(\mathcal{G}_v)\simeq \piet(X^v_0)$.
				\end{minipage}
				
				\smallskip
				
				This construction is a generalization of Serre's tree group theory in the sense that (CombAG) doesn't require fixing a compatible set of base points for the various local groups involved.

			\subsubsection{A combinatorial Grothendieck conjecture-type}\label{sub:CombAG} 
				A typical example of combinatorial Grothendieck conjecture-type result is given as follows -- see \cite{CbTpII} Th.~1.9 for an exact statement.
				\begin{theo}\label{th:CombGC} Consider $\mathcal{G}$ and $\mathcal{H}$ two semi-graphs of anabelioids of PSC-type, equipped with outer representations $\rho_{\mathcal{G}}\colon I\to \Out(\Pi_{\mathcal{G}})$ and $\rho_{\mathcal{H}}\colon I\to \Out(\Pi_{\mathcal{H}})$ as above\textsuperscript{$\dagger$}, where $I$ denotes an inertia group. Then any isomorphism $\alpha\colon \Pi_{\mathcal{G}}\simeq  \Pi_{\mathcal{H}}$ compatible with $\rho_{\mathcal{G}}$ and $\rho_{\mathcal{H}}$ determines an isomorphism of semi-graphs of anabelioids $\tilde{\alpha}\colon \mathcal{G}\simeq \mathcal{H}$ -- and, in particular, also between the components at edges, vertices, and nodes.
				\end{theo}
				\hspace{-1em}\hphantom{A}\textsuperscript{$\dagger$}More exactly, one must further assume that the outer representations $\rho_{\mathcal{G}}$ and $\rho_{\mathcal{H}}$ satisfy certain group-theoretic properties on the inertia groups at nodes -- that they are ``Nodally Non-degenerate'' (NN), or of ``Inertia PSC-type'' (IPSC). For an outer representation to be of IPSC-type essentially means that it is of the type of $\widehat{\mathbb{\ZZ}}$-inertia action that comes from the Galois-monodromy of a curve $X/\mathcal{O}_k$, that is -- see ibid.:
				\[
				\begin{tikzcd}
					I \ar[d,"\sim"',sloped]\ar[r]& \Out(\Pi_{\mathcal{G}})\ar[d,"\sim",sloped]\\
					G_{\bar{\eta}/\eta} \ar[r]\ar[phantom,"\circlearrowleft",ru]& \Out(\Pi_{\mathcal{G}_{X_s}})
				\end{tikzcd}
				\]
				Note that $\alpha$ is only required to be compatible with $\rho_{\mathcal{G}}$ and $\rho_{\mathcal{H}}$ up to an automorphism of $I$, so that the statement is of \emph{semi-absolute} rather than relative nature.
				
				\medskip
				
				(CombAG) deals with the group-theoretic characterization of anabelian data in higher dimensional varieties, such as configuration spaces of hyperbolic curves, as presented in \S~\ref{sec:GTAnab} and \S~\ref{sec:BGT}. It further allows one to establish properties of \emph{profinite} Dehn twists and their conjugacy classes \cite{CbTpI} \S~4-5, and provides certain type of canonicity with respect to tangential structures (synchronization of cyclotomes, that is of groups isomorphic to $\widehat{\ZZ}$, see ibid. \S~3).
				
				\medskip
				
				We refer to the survey \cite{Hos22} \S~7-\S~8 for a recent report on \emph{Applications of and Progress in Combinatorial Anabelian Geometry}.

	\section{Anabelian arithmetic of the genus zero modular tower}
		We now consider the  \emph{Teichmüller tower} $\mathcal{T}=\{\piet(\MM_{0,m}\times\bar{\QQ})\}_{m\geq 3}$, whose morphisms are induced by the Knudsen point-retraction $\MM_{0,m+1}\to\MM_{0,m}$. We further denote by $\Out_\mathcal{S}(\mathcal{T})$ the tower of outer automorphisms that \textbf{(a)}~commute with the point-retraction morphisms and, in order to reflect the full geometry of $\MM_{0,[m]}$ as in \S~\ref{sub:Hsym}, \textbf{(b)}~that level-wise commute with the $S_m$-action on points.
		
		\medskip
		
		Since the previous setup is $\QQ$-compatible, one obtains a morphism $G_\QQ\to \Out_\mathcal{S}(\mathcal{T})$, which is already injective at the level $m=4$ by Belyi. In this context, Grothendieck's original insight of \cite{GroEsq97} regarding the surjectivity of this morphism, can be formulated, with $\GT\geq G_\QQ$ as in \S~\ref{sub:Hsym}, as combinatorial or group-theoretic container, as the question: 
		
		\medskip
		
		{\centering\itshape 
			``To what extent does one have an identification  $\GT\simeq \Out_\mathcal{S}(\mathcal{T})$?''
		\par} 
		
		\medskip
		
		A first partial answer is provided by the introduction of a subgroup $\Out_\mathcal{S}^\sharp(\mathcal{T})\leq \Out_\mathcal{S}(\mathcal{T})$ consisting of \emph{quasi-special automorphisms}, that satisfy a certain additional arithmetic property -- see Th.~\ref{th:GTHS5} and compare to the identification of Eq.~\eqref{eq:GTanab}.
		
		\medskip
		
		The naturality of this quasi-special arithmetic property stems from anabelian considerations concerning $\mathcal{M}_{0,m}$ as in \cite{Nak94}, for which the Grothendieck–Teichmüller group arises as a potentially ideal group-theoretic container.
		
		\subsection{Anabelian structures on $\MM_{0,m}$} 
			In the Deligne-Mumford compactification of $\MM_{0,m}$, let $D_{ij}$, $1\leq i<j\leq m$, denote the irreducible divisor of the boundary that supports stable $m$-pointed curves of genus $0$ with one singular point and the marked points $x_i$ and $x_j$ on a same irreducible component. 
			
			\medskip
			
			In the geometric étale fundamental group $\piet(\MM_{0,[m]}\times \bar{\QQ})\simeq \widehat{\Gamma}_{0,[m]}$, a loop around $D_{ij}$ corresponds to a pure braid generator $x_{ij}\in \widehat{P}_m$, while, at the arithmetic level, the absolute Galois group $G_\QQ$ conjugates and permutes the associated cuspidal inertia groups $\langle x_{ij}\rangle$.
			
			\medskip
			
			We thus denote by $\mathfrak{X}_{ij}$ the conjugacy class of such cuspidal inertia subgroups as in \cite{Nak94} \S~3.1.
			
			\subsubsection{Cuspidal inertia and quasi-special automorphisms} 
				At the group-theoretic level, the outer automorphisms of $\widehat{\Gamma}_{0,m}$ that reflect the above Galois property -- that is that fix the $\mathfrak{X}_{ij}$ -- are denoted $\Out^\flat(\widehat{\Gamma}_{0,m})$ and called \emph{quasi-special}.
				
				\medskip
				
				Since the $\{ \mathfrak{X}_{ij}\}_{1\leq i\leq m}$ generates the kernels of the point-contraction induced morphisms $p_i\colon \widehat{\Gamma}_{0,m}\to \widehat{\Gamma}_{0,m-1}$, one obtains a series of canonical homomorphisms
				\[
					p_i\colon\Out^\flat(\widehat{\Gamma}_{0,m})\to \Out^\flat(\widehat{\Gamma}_{0,m-1})\quad \text{ for }1\leq i\leq m,
				\]
				so that, testing equality on two components and by intricate pure braid computations with respect to a well-chosen presentation of $P_m$ as in ibid. 3.1.3.1-7, one obtains, see ibid. Lem.~3.2.2.
				
				\begin{prop}\label{prop:NakOutInj}
					The morphisms
					 \[
					 	(p_i,p_j)\colon\Out^\flat(\widehat{\Gamma}_{0,m})\to \Out^\flat(\widehat{\Gamma}_{0,m-1})\times \Out^\flat(\widehat{\Gamma}_{0,m-1})
					 \]
					  are injective for $1\leq i,j\leq m$, $i\neq j$.
				\end{prop}
				We refer also to Prop.~\ref{prop:HSOutinj} and Th.~\ref{theo:OutFTinj} for similar injective results.
				
				\medskip
				
				This property allows to reduce establishing certain arithmetic geometric statements to the case $m=4$ -- that is to $\MM_{0,4}\simeq \mathbb{P}^1\setminus\{0,1,\infty\}$ -- such as for Th.~\ref{th:M0manabEll} below.
			
			\subsubsection{Anabelian properties}\label{subsub:AnabNak} 
				Let us show, following  the original \cite{Nak94} for the pro-$\ell$ completion -- see Theorem~A ibid. -- how the previous group-theoretic input establishes \emph{(GC) for the higher-dimensional $\MM_{0,m}$s} in its full profinite version, that is:
				\begin{theo}\label{th:M0manabEll} 
					For $m\geq 4$  the natural application
					\begin{equation}\label{eq:GCM0m}
						\psi\colon\Aut_\QQ(\MM_{0,m})\longrightarrow \Aut_{G_\QQ}(\piet (\MM_{0,m}))/\operatorname{Inn} \piet(\MM_{0,m}\times \bar{\QQ}).
					\end{equation}
					is a bijection.
				\end{theo}

			This result follows from some classical but elaborate group-theoretic arguments such as the group-theoretic characterization of cuspidal inertia groups of a hyperbolic curve $X$ as ``cyclotomic normalized'' subgroup of $\piet(X\times \bar{K})$ -- see \cite{Nak94} Th.~2.1.1, which is a refined version of the Deligne's weight characterization of cuspidal groups developed in \cite{Nak90}: profinite group theory for surface groups associates to the homotopy (FES) of $X$ a field of definition $k_0$, which in turns provides: 
			\[
				\text{A $\ell$-adic representation } G_{k_0}\to \operatorname{GL}(H_1(X_{\bar{k}_0},\mathbb{Z}_{\ell}))\text{ where {\itshape the cuspidal groups have weight}} -2 
			\]
			see also \cite{Nak90} \S~(2.3). Note that this technique is still today part of the anabelian geometer toolbox, see \cite{CPY26} as one of many examples.
			
			\medskip
			
			In the case of $\MM_{0,m}$, this provides, with the use of the dimension-shifting morphisms $p_i$, a group-theoretic retraction of $\psi$ (that is $\widetilde{\psi}\circ\psi=\operatorname{id}$),
			\[
			\widetilde{\psi}\colon \Aut_{G_\QQ}(\piet(\MM_{0,m}))/\operatorname{Inn} \piet(\MM_{0,m}\times \bar{\QQ})  
				\longrightarrow
			\Aut_\QQ(\MM_{0,m}),
			\]
			which to a $G_\QQ$-equivariant automorphism $\phi$ of $\piet(\MM_{0,m})$ attaches a $\QQ$-automorphism of $\MM_{0,m}$ whose permutation on the marked points is identical to the one induced by $\psi$ on the cuspidal subgroups -- see Th.~3.1.13 ibid., where the same series of arguments holds for the profinite completion, since also a full class of finite groups.      
      
      \medskip

      The injectivity of $\widetilde{\psi}$ -- equivalently the surjectivity of $\psi$ -- as well as the fact that the bijectivity can thus be reduced to the case $m=4$, follows from Prop.~\ref{prop:NakOutInj} as applied to the centralizer $Z_{\Out(\widehat{\Gamma}_{0,m})}(G_\QQ)$. 
      
      The bijectivity for $m=4$ is then given by the (GC) of \cite{Tam97}, which is a strongest form of the anabelian result given in~\cite{Nak91}.
		
		\subsection{The group $\GT$ as a group-theoretic container}\label{sub:GTtower}
			In \cite{HS00}, the Grothendieck-Teichmüller group is used as a container to better seize the anabelian properties of $\piet(\MM_{0,m}\times\bar{\QQ})$ of the previous section at a group-theoretic level. 
			
			\medskip
			
			For example, the ``nearly anabelian'' property that $\Aut_\QQ(\MM_{0,m})\simeq S_m$ for $m\geq 5$ (see also below Eq.~\eqref{eq:GTanab}) -- or even the construction of the section of Eq.~\eqref{eq:GCM0m} -- and the braid group property that the $\GT$-action restricts from $\widehat{B}_m$ to $\widehat{P}_m$ (with $\widehat{B}_m/\widehat{P}_m\simeq S_m $), see discussion below Th.~\ref{th:GTBn}, motivates to consider the subgroup $\OutNS{m}\leq\Out^\flat(\widehat{\Gamma}_{0,m})$ of outer automorphisms that commute with the $S_m$-action on points -- the so-called \emph{quasi-special and symmetric} automorphisms.

			\begin{wrapfigure}[6]{r}{5cm}
				\vspace*{-1em}
				\centering
				$\begin{tikzcd}[column sep=1em]
					G_\QQ \ar[r,hook] \ar[rd]& \Out^\sharp(\widehat{\Gamma}_{0,5})\simeq \GT\ar[d,dashed,shift left=.9ex,"s_5", bend left]\\
					& {\rm A}_5\leq \Aut(\widehat{\Gamma}_{0,5})\ar[u,twoheadrightarrow,shift left=.9ex]
				\end{tikzcd}$
				\captionof{figure}{Belyi lifting: $s_5$ is a section of $\mathrm{Aut}\to\mathrm{Out}$}
			\end{wrapfigure}
			\subsubsection{Belyi lifting and the identification at level $5$} 
				Considering further the \emph{Belyi lifting $A_5$} -- that is defined as the group of automorphisms of elements that acts similarly to the Ihara $G_\QQ$-action on $\widehat{\Gamma}_{0,[5]}$ when considered as the geometric fundamental groups $\piet(\MM_{0,5}\times\bar{\QQ})$ -- some braid group theoretic arguments similar to \cite{Nak94} as in \S~above, establish the existence of a section $s_5$ -- see \cite{HS00} \S~2. 
			
				\medskip
				
				With additional computation in terms of the cocycle definition of $\GT$ as in Rem.~\ref{rem:GTCohom} one identifies as in ibid. Th.~3.
				\begin{theo}\label{th:GTHS5} 
					The Belyi lifting $A_5$ and the associated section $s_5\colon \OutNS{5}\to A_5$ induce an isomorphism
					\begin{equation}
            \psi\colon\Out^\sharp(\widehat{\Gamma}_{0,5})\simeq\GT
					\end{equation} 
				\end{theo}

				This result provides another precise formulation of the insight that \emph{the group $\GT$ -- and thus the action of $G_\QQ\hookrightarrow \GT$ on the modular tower $\mathcal{T}$ -- is defined via the first two levels of the tower only}.

			\subsubsection{Higher automorphisms}\label{subsec:Haut} 
				As in the previous section, the quasi-special property of automorphisms tantamounts to preserving the homotopy kernels of the Knudsen-Mumford contraction morphisms, so that one obtains -- see ibid. Prop.~8.
				
				\begin{prop}\label{prop:HSOutinj}
					For $m\geq 5$, there exist natural homomorphisms
					\begin{equation}\label{eq:OutInjHS}
						q_i\colon\OutNS{m}\longrightarrow\OutNS{m-1} \text{ for } i=1,\dots,m
					\end{equation}	
					and these homomorphisms are injective.
				\end{prop}
				
				\medskip
				
				This result follows more precisely from the Birman point-erasing exact sequence, which algebraically splits on generators, so that $\widehat{\Gamma}_{0,m}\simeq \widehat{\mathbb{F}}_2\rtimes\widehat{\Gamma}_{m-1}$ fits in a commutative
				
				\medskip
				
				{\centering
					\begin{tikzcd}
						1 \ar[r]& \widehat{\Gamma}_{0,m} \ar[r]\ar[d]& \Aut_m^\sharp \ar[r]\ar[d]& \OutNS{m} \ar[r]\ar[d,"q_i"]& 1\\
						1 \ar[r]& \widehat{\Gamma}_{0,m-1} \ar[r]& \Aut_{m-1}^\sharp \ar[r]& \OutNS{m-1} \ar[r]& 1
					\end{tikzcd}
				\par}
				
				\medskip
				
				The injectivity follows from some reasoning on the permutation $S_m\to \Out(\widehat{\Gamma}_{0,m})$ and from the injectivity result of Prop.~\ref{prop:NakOutInj}.
			
			\begin{wrapfigure}[10]{r}{5cm}
				\centering
				\vspace{-.5em}
				$\begin{tikzcd}[row sep=1.2em, /tikz/column 1/.style={column sep=1em}]
					\AutNS{m} \ar[d]\ar[r]&  \OutNS{m} \ar[d,hook,"q_{m}"']& \\
					\AutNS{m-1} \ar[d,dotted]\ar[r]& \OutNS{m-1} \ar[d,dotted,hook,"q_{m-1}"']& \GT \ar[dotted,hook,lu, "e_m"']\ar[dotted,hook,l, "e_{m-1}"'] \ar[dotted,hook,ld, "e_5"']\\
					\AutNS{5} \ar[d]\ar[r]& \OutNS{5} \ar[d,hook,"q_{5}"]& \\
					\AutNS{4} \ar[r]& \OutNS{4}& 
				\end{tikzcd}$
				\captionof{figure}{$\GT$ and tower}\label{fig:GTTower}
			\end{wrapfigure}
			\subsubsection{Tower morphisms and outer stability}
				The stability of the outer-tower, that is a group-theoretic result, is then obtained as a consequence of the group-theoretic Galois isomorphism of Th.~\ref{th:GTHS5} --- see Th.~5 ibid.
				
				\begin{theo}\label{theo:HSMain}
					For $m\geq 5$, the map
					
					\medskip
					
					{\centering
						$q_m\colon\OutNS{m}\to \OutNS{m-1}$
						\par}
					
					\medskip
						
					is an isomorphism. As a consequence, $\GT\simeq \OutNS{m}$ for $m\geq 5$.
				\end{theo}

				One checks immediately that the classical $\GT$ action on braids of Th.~\ref{th:GTBn} induces homomorphisms 
				
				\medskip
				
				{\centering
					$e_m\colon \GT \longrightarrow \OutNS{m}$, for $m\geq 4$
					\par}
				
				\medskip
				
				that are each compatible with the morphism $q_m$. 
				
				\medskip
				
				In particular, the compatibility at level $m=5$ implies that \emph{the morphism $e_5$ is an inverse isomorphism of $\psi$} in Th.~\ref{th:GTHS5}, which in turns implies that
				
				\medskip
				
				{\centering \itshape
					Every homomorphism $e_n\colon\GT\to \OutNS{n}$ is an isomorphism.
				\par}
							
			\begin{rmk}\label{rem:GTCohom}
				In a manner that brings $\GT$ closer to group-theoretic anabelian methods, one observes that the defining relations (I), (II), and (III) of $\GT$ may be interpreted as cocycle relations in the corresponding non-commutative cohomology sets $H^1(\langle \theta \rangle,\widehat{\mathbb{F}}_2)$, $H^1(\langle \omega\rangle,\widehat{\mathbb{F}}_2)$, and $H^1(\langle \rho\rangle,\widehat{\Gamma}_{0,5})$ respectively -- see \cite{LS97}.
			\end{rmk}
			
			\bigskip
			
			While this Grothendieck-Teichmüller approach allows one to exploit the $S_m$-anabelian properties -- see below Eq.~\eqref{eq:GTanab} for the anabelian significance of $S_m$ -- in a more group-theoretically consistent manner, the results still essentially derive from anabelian geometry, and the techniques remain deeply anchored in \emph{ad hoc} braid-theoretic computations. At this stage, the precise nature of the $\flat/\sharp$ condition -- arithmetic or genuinely group-theoretic -- remains open.

	\section{The group $\GT$ as anabelian container}\label{sec:GTAnab}
		We reach what constitutes a double breakthrough, that is the development of the purely combinatorial anabelian techniques of \S~\ref{sub:whatisCbAn}-\ref{sub:CombAG} for \textbf{(a)}~the identification of $\GT$ as an anabelian object and \textbf{(b)}~the removal of the $\flat/\sharp$ condition (thus confirming their purely group-theoretic nature) -- see \cite{HMM22}.
		
		\begin{theo}
			The Grothendieck-Teichmüller group is anabelian, in the sense that it identifies with anabelian objects
			\[
				\GT\simeq Z_{\Out(\piet(\MM_{0,5}\times \bar{\QQ}))}[\varinjlim \{Z_{\Out(\piet(\MM_{0,5}\times \bar{\QQ}))}(H),\, H\leq_{\text{open}}\Out(\piet(\MM_{0,5}\times \bar{\QQ}))\}],
			\]
			where the direct limit is the so-called local center $Z^{\rm loc}(\Out(\piet(\MM_{0,5}\times \bar{\QQ})))$, which can be identified with $S_5$.
		\end{theo}
		Writing $\Pi_m=\piet(\conf{m})$, this result is indeed a corollary of the isomorphism below, which is given as Cor.~2.8 ibid.:
		\begin{equation}\label{eq:GTanab}
			\widehat{GT}\times S_{m+3}\simeq \Out(\Pi_m)\simeq \Out(\widehat{\Gamma}_{0,m+3})\text{ for } m\geq 2,
		\end{equation}
    which also establishes, unconditionally, that $\GT$ and $S_{m+3}$ are mutual centralizers in $\Out(\Pi_m)$. The Galois centralizer is moreover $Z_{\Out(\piet(\MM_{0,m}\times \bar{\QQ}))}(G_\QQ)\simeq S_m$ by Th.~\ref{th:M0manabEll}.
		
		\medskip
		
\begin{wrapfigure}[9]{r}{6cm}
					\centering
					\vspace{-1.5em}
					\def\svgwidth{\linewidth}
					\resizebox{.7\linewidth}{!}{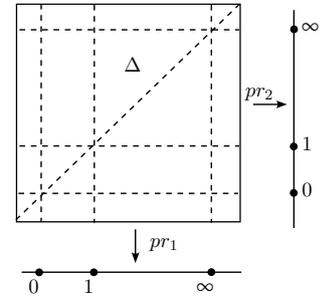}
					\captionof{figure}{Tripods synchronization in $X_2$}\label{fig:TripSync}
				\end{wrapfigure}

		In the following, we identify the genus zero moduli spaces of curves $\MM_{0,m+3}\times\bar{\QQ}\simeq (\mathbb{P}^1_{\bar{\QQ}}\setminus\{0,1,\infty\})^{m}-\Delta$, that is with the configuration space $\conf{m}$ of the geometric projective line $X_{0,3}=\mathbb{P}^1_{\bar{\QQ}}\setminus\{0,1,\infty\}$. Through this section, the reader may keep in mind the first non-trivial case $m=2$, where $\Pi_2=\piet(\operatorname{Conf}_2(X_{0,3}))$ and $\MM_{0,5}\simeq\operatorname{Conf}_2(X_{0,3})$ is the surface drawn in Fig.~\ref{fig:TripSync} -- see \S~\ref{sub:TripSync}.
		
		\subsection{Decuspidalization of automorphisms}
			The notion of \emph{FC-admissible automorphisms} introduced in \cite{CmbCsp} Def.~1.1~(ii) is equivalent to the one of quasi-special automorphisms, that is $\Out^\flat_m=\Out^{\rm FC}_m$ (see ibid. Prop. 1.3~(vii)),  but formulated in a way that is closer to the anabelian mechanisms.

			\subsubsection{Cusps in Fiber-admissible automorphisms}
			The label-forgetting projections $X_n\to X_m$ define a subgroup filtration $K_m=\ker(\Pi_n\twoheadrightarrow \Pi_m)\geq K_{m+1}$ which fits in the diagram
			\[
			\begin{tikzcd}
				1 \ar[r]& K_m/K_{m+1} \ar[r]\ar[d,"\sim",sloped]& \Pi_n/K_{m+1} \ar[r] \ar[r]\ar[d,"\sim",sloped]& \Pi_n/K_{m} \ar[r] \ar[r]\ar[d,"\sim",sloped]& 1\\
				1\ar[r] &  \pi_{0,m+3} \ar[r] & \Pi_{m+1}\ar[r] &\Pi_{m} \ar[r] & 1
			\end{tikzcd}
			\]
			where $\pi_{0,m+3}$ is the geometric fundamental group of a hyperbolic curve of type $(0,m+3)$, so that it contains in particular some cuspidal subgroups. The bottom exact sequence can be seen as a refined version of the Birman point-erasing sequence used in \S~\ref{subsec:Haut}.
			
			\medskip
			
			Accordingly, the group $\Out^{\rm FC}(\Pi_m)<\Out^{\rm F}(\Pi_m)<\Out(\Pi_m)$ of \emph{FC-admissible automorphisms}\footnote{This notion, and most of the results of this section -- e.g. Th.~\ref{th:FCF} with $m\geq 3$ for $r>0$, and $m\geq 4$ for $ r=0$ -- extend readily to configuration spaces of curves of any hyperbolic type $(g,r)$, that is $2g-2+r> 0$, where $\Pi_m$ denotes the fundamental group of configuration space of $m$ points on curves of type $(g,r)$.} is composed of outer automorphisms that preserve the kernels $K_m$ (Fiber-admissible), and that induce a permutation of the cuspidal subgroups of $\pi_{0,m+3} \simeq K_m/K_{m+1}$ (Cuspidal admissible).

			\subsubsection{The FC-F automorphisms identification}
			While the Fiber-admissible property ensures the existence of morphism at the outer-level, a major result of \cite{CbTpII} Th.~A and Th.~2.3~(ii) is that the Cuspidal condition is automatically satisfied.
			\begin{theo}\label{th:FCF}
				For $m\geq 2$, we have $\Out^{\rm FC}(\Pi_m)\simeq \Out^{\rm F}(\Pi_m)$.
			\end{theo}
		
			Rather than a sketch of the proof of this result, we present the reader with a second result, also of independent arithmetic interest, which relies on identical arguments and which reveals a broader scope of the combinatorial anabelian geometry techniques -- see Th.~\ref{theo:OutFTinj}.
			
			\subsection{Generalized Belyi injectivity} \label{sub:BelInj}
				Beyond its intrinsic interest, the following result implies the generalization, from Belyi's original injectivity for $X=\mathbb{P}^1\setminus\{0,1,\infty\}$ to the injectivity for any hyperbolic curve $X$ of type $(g,m)$ of
				\[
					G_K\longrightarrow \Out(\piet(X\times\bar{K})),
				\]
				see Th.~\ref{theo:BelyiP} below.
			
				\medskip
				
				The following is \cite{CbTpII} Th.~A and Th.~2.3~(i) and (ii).
				\begin{theo}\label{theo:OutFTinj}
					For $X$ a hyperbolic curve over an algebraically close field of characteristic zero, the projection morphisms 
					\[
						p_i\colon\Out^{\rm FC} (\Pi_m)\longrightarrow \Out^{\rm FC} (\Pi_{m-1})\text{ for } i=1,\dots,m
						\]
					are injective for any  $m\geq 2$, and bijective for $m\geq 4$.
				\end{theo}

				\medskip
				
				We now give a sketch of the involved mechanisms.
				
				\subsubsection{Tripod synchronization and (CombGC)}\label{sub:TripSync} Let us fix $m=2$ for simplicity. Consider $\alpha\in \Out^{\rm FC}(\Pi_2)$ which by assumption becomes, up to inner automorphism, the identity on $\Pi_1$ by one of the surjections. The three key steps are as follows.
				
				\medskip
				
				\noindent(i)~\emph{Synchronization of tripods.} A tripod is a divisor that is isomorphic to $\mathbb{P}^1\setminus\{0,1,\infty\}$. Let $I_1$, $I_2$ and $I_\Delta$ be the tripod subquotients of $\Pi_2$ attached to the horizontal, the vertical, and the diagonal tripod divisors of $\MM_{0,5}\simeq \conf{2}$ (three out of the seven), each isomorphic to the fundamental group of a tripod. A first result of (CmbAG) is to establish the outer-synchronization of these tripods
				
				\medskip
				
				{\centering
				
				$\Out(I_1)\xrightarrow{\sim} \Out(I_\Delta) \xleftarrow{\sim} \Out(I_2)$,
				\par
				}
				
				\medskip
				
				under which $\alpha|I_1$ and $\alpha|I_2$ both correspond to $\alpha|I_\Delta$, and which implies that $\Pi_{1/2}\coloneq\operatorname{Ker}(p_1)\simeq \operatorname{Ker}(p_2)$ is independent of the projection.
				
				\medskip
				
				\noindent(ii)~\emph{Application of (CombGC).} Consider the fiber $F_0$ at $0\in X$ of any of the projections $pr_i\colon X_2^{\rm log}\to X^{\rm log}$, where the schemes are endowed with the log-structure associated to the normal crossing divisors of their smooth compactification. 
				\medskip
				
				Notice that $\Pi_{1/2}\simeq\pi_1^{\rm log}(T\cup F_0)$ where $T\cup F_0$ is a stable curve with a tripod exceptional divisor $T$, with $\Pi_{F_0}\leq \Pi_{2/1}$, so that one has the series of commutative squares:
				
				\[
				\begin{tikzcd}
					1 \ar[r] & \{\piet(F_0),\piet(T)\} \subset \ker\psi \ar[r]\ar[d,"\simeq",sloped, phantom]&\pi_1^{\rm log}(T\cup F_0) \ar[r, "\psi"] \ar[d,hook]&\pi_1^{\rm log}(\{0\}) \ar[r] \ar[d,hook]& 1\\
					1 \ar[r]& \Pi_{2/1} \ar[r]& \Pi_2 \ar[r] & \Pi_1 \ar[r] & 1
				\end{tikzcd}
				\]
				
				\begin{wrapfigure}[6]{r}{6cm}
					\centering
					\vspace{-1em}
					\def\svgwidth{\linewidth}
					\resizebox{.7\linewidth}{!}{
\begingroup%
  \makeatletter%
  \providecommand\color[2][]{%
    \errmessage{(Inkscape) Color is used for the text in Inkscape, but the package 'color.sty' is not loaded}%
    \renewcommand\color[2][]{}%
  }%
  \providecommand\transparent[1]{%
    \errmessage{(Inkscape) Transparency is used (non-zero) for the text in Inkscape, but the package 'transparent.sty' is not loaded}%
    \renewcommand\transparent[1]{}%
  }%
  \providecommand\rotatebox[2]{#2}%
  \newcommand*\fsize{\dimexpr\f@size pt\relax}%
  \newcommand*\lineheight[1]{\fontsize{\fsize}{#1\fsize}\selectfont}%
  \ifx\svgwidth\undefined%
    \setlength{\unitlength}{66.13788455bp}%
    \ifx\svgscale\undefined%
      \relax%
    \else%
      \setlength{\unitlength}{\unitlength * \real{\svgscale}}%
    \fi%
  \else%
    \setlength{\unitlength}{\svgwidth}%
  \fi%
  \global\let\svgwidth\undefined%
  \global\let\svgscale\undefined%
  \makeatother%
  \begin{picture}(1,0.54015009)%
    \lineheight{1}%
    \setlength\tabcolsep{0pt}%
    \put(0,0){\includegraphics[width=\unitlength,page=1]{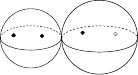}}%
    \put(0.06141327,0.33142498){\color[rgb]{0,0,0}\makebox(0,0)[lt]{\smash{\begin{tabular}[t]{l}$\infty$\end{tabular}}}}%
    \put(0.58097256,0.36113557){\color[rgb]{0,0,0}\makebox(0,0)[lt]{\smash{\begin{tabular}[t]{l}$0$\end{tabular}}}}%
    \put(0.28607206,0.33561621){\color[rgb]{0,0,0}\makebox(0,0)[lt]{\smash{\begin{tabular}[t]{l}$1$\end{tabular}}}}%
  \end{picture}%
\endgroup%
}
					\captionof{figure}{Fiber of $X_2\to X_{0,3}$}
				\end{wrapfigure}
				Since $\Pi_{F_0}$ contains an inertia group $I_0\simeq \widehat{\mathbb{\ZZ}}$ and since $\alpha_{2/1}$ is of IPSC-type, Th.~\ref{th:CombGC} implies that $\alpha_{2/1}$ stabilizes $\Pi_{F_0}\simeq \Pi_1$. A similar argument for the fiber $F_1$ implies that
				
				\medskip
				
				{\centering
					$\alpha_{2/1}|\Pi_{F_0}=\alpha_{2/1}|\Pi_{F_1}={\rm id}$
				\par	
				}
					
				\medskip
					
				(iii)~\emph{Van Kampen.} Since $\Pi_{2/1}$ is obtained by gluing over the various fibers, we conclude the generic triviality $\alpha_{2/1}={\rm id}$ from the identity at each of the fibers.
				
				\begin{rmk}
					A similar sequence of arguments with similar conclusion holds for step (ii) in the case where $X$ is hyperbolic of type $(g,r)$ with $r\geq 2$.
				\end{rmk}
				
				\medskip
				
				\begin{wrapfigure}[10]{r}{6cm}
					\centering
					\vspace{-2.5em}
					$\begin{tikzcd}[/tikz/column 2/.style={column sep=.5em}]
						G_K \ar[r, "\alpha_N"]\ar[rd, "\alpha_{N-1}"{pos=0.75},rounded corners, to path={[xshift=10ex]|- (\tikztotarget.west)\tikztonodes}]\ar[rdd,"\alpha_4"{pos=0.75}, rounded corners, to path={[xshift=10ex]|- (\tikztotarget.west)\tikztonodes}] \ar[rddd,"\alpha_3"{pos=0.8}, rounded corners, to path={|- (\tikztotarget.west)\tikztonodes}] &  \Out^{\rm FC}(\Pi_N) \ar[d,hook]& \\
						& \Out^{\rm FC}(\Pi_{N-1}) \ar[d,dotted,hook]& \\
						& \Out^{\rm FC}(\Pi_4) \ar[d,hook]& \\
						& \Out^{\rm FC}(\Pi_3)\ar[r,"T"]&\Out(\piet(X_{0,3}))
					\end{tikzcd}$
					\captionof{figure}{Generalized Belyi injectivity}
				\end{wrapfigure}
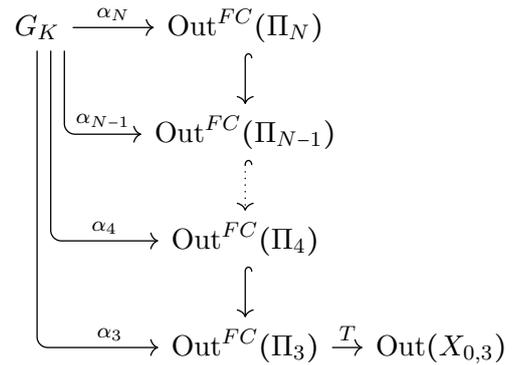
				\subsubsection{The injectivity result} The following is the abutment of a long series of results, from Belyi injectivity for genus zero hyperbolic curves (1979) to Hoshi-Mochizuki \cite{NodNon11} Cor.~6.3, via Voevodski in genus one (1991), Matsumoto for affine curves (1996) -- see references in ibid.:
				\begin{theo}\label{theo:BelyiP}
					Let $X$ be a hyperbolic curve of type $(g,m)$ over a number field $K$, then the morphism $ G_K\to \Out(\Pi_N)$ is injective.
				\end{theo}
				One shows that, by restriction to $\Out^{\rm FC}(\Pi_N)\leq\Out(\Pi_N)$. This follows from Th.~\ref{theo:OutFTinj}, that is from the injectivity of 
				
				\medskip
				
				{\centering
					$\Out^{\rm FC}(\Pi_N)\to \Out^{\rm FC}(\Pi_{N-1})$ 
				\par}
				
				\medskip
				
				which reduces the results to the case $N=3$ and thus to Belyi's original one -- since $\alpha_3\circ T$ is Belyi's original morphism to $X_{0,3}=\mathbb{P}_{\bar{\QQ}}^1\setminus\{0,1,\infty\}$. Here $T$ denotes the \emph{tripod homomorphism} of \cite{CbTpII} Def.~3.19, which attaches to an FC-admissible outer automorphism the one it induces on the \emph{central tripod} -- for $N=3$, the tripod obtained by blowing up the intersection of the three diagonals of $X^3$ -- and whose well-definedness follows the tripod synchronization above.
				
			\subsection{The group $\GT$ is anabelian}
				In order to exhibit the anabelian nature\footnote{We let the attentive reader to write the group-theoretic algorithm that defines $\GT$ from $\Pi_m$.} of $\GT$, we rely on the group-theoretic notion of automorphisms which respect the forgetting of marked points: an automorphism of $\Pi_m$ is \emph{gF-admissible}  (gF for ``generalized fiber'') if it stabilizes any $\ker(\Pi_m\to \Pi_n)<\Pi_m$ obtained by forgetting $m-n$ points; we denote by $\Out^{\rm gF}(\Pi_m)<\Out^{\rm F}(\Pi_m)$ the corresponding outer automorphism group. 
				
				\medskip
				
				One then shows that such automorphisms permute the generalized fiber subgroups of co-length $m-n=1$, so that one obtains as in \cite{HMM22} Cor.~2.6, an exact sequence:
				\[\begin{tikzcd}
					1\ar[r, dotted]& \Out^{\rm gF}(\Pi_m)\ar[r, dotted]  &\Out(\Pi_m)\ar[r, "\exists"]& S_{m+3}\ar[r] & 1
				\end{tikzcd}	
				\]
				with a natural splitting given by the permutation of points of $\MM_{0,m+3}$, so that:
				\begin{equation}\label{eq:gfDir}
					\Out^{\rm gF}(\Pi_m)\times S_{m+3} \simeq \Out(\Pi_m)
				\end{equation}
				
				\medskip
				
				The final $\GT$-isomorphism of Eq.~\eqref{eq:GTanab} then results as follows:
				\begin{align*}
					\GT &\simeq\Out^{\rm FC}(\Pi_m)\cap Z_{\Out(\Pi_m)}(S_{m+3}) && \text{by rewriting Th.~\ref{theo:HSMain},}\\
					  &\simeq\Out^{\rm F}(\Pi_m)\cap \Out^{\rm gF}(\Pi_m)  && \text{since $Z(S_{m+3})=\{1\}$ in Eq.~\eqref{eq:gfDir} and by Th.~\ref{th:FCF},}\\
					\GT &\simeq\Out^{\rm gF}(\Pi_m) && \text{so that one applies Eq.\eqref{eq:gfDir} again to conclude.}
				\end{align*}

        \medskip

        In particular, $\GT$ appears as an open normal subgroup of $\Out(\Pi_m)$, of index $(m+3)!$.

				\medskip
				
				\begin{rmk}\label{rem:grpAdic}
					Results of this section still rely on some explicit $\ell$-adic Lie algebra computations via the group presentations of \cite{Nak94} in the following way:
					\begin{enumerate}
						\item  The stabilization of generalized fiber subgroups relies on the isomorphism
						\[
						\Aut(\operatorname{Gr}(\Pi_m^\ell))\simeq S_{m+3}\times \QQ_\ell^\times
						\]
						where the grading is associated to the lower central series; 
						\item The trivial center property in the final sequence of isomorphisms follows from Th.~A ibid.
					\end{enumerate}
					The group-theoretic reconstruction of the discrete invariants, that are the genus $g$, the number of cusps $m$ and the dimension $N$ of configuration spaces of $X_N$ also follows some $\ell$-adic Lie algebra computations -- see \cite{HMM22} Theorem~2.5, and even more generally \cite{Saw16}.
				\end{rmk}

	\section{Combinatorial models of $G_\QQ$}\label{sec:BGT}				
		We now sketch the anabelian algorithm that reconstructs $\GT$ from $\Pi_m$, which, in a certain sense, contracts the arithmetic-geometry-group theory domain delineated by the triangle of Fig.~\ref{fig:AGtr} -- In the following we write $X_{0,3}=\mathbb{P}_{\bar{\QQ}}^1\setminus\{0,1,\infty\}$ and  $X_m=\operatorname{Conf}_m(X_{0,3})$.
		
		\begin{wrapfigure}[8]{r}{5cm}\vspace{-.5em}
			\centering\small
			\begin{tikzcd}[row sep=4pt, column sep=15pt,/tikz/row 2/.append style={row sep=20pt}, ampersand replacement=\&]
				V \ar[ddr, twoheadrightarrow, bend right, "f.et."'] \ar[dr,twoheadrightarrow, dotted,  "fin.","et."']\& \& U_X\ar[ddl,hook, dotted, bend left,"op."] \ar[dl,twoheadrightarrow, dotted,  "\exists\ fin."', "et."]\\
				\&\mathbb{P}^1\setminus D  \& \\
				\& X \&	\\[-7pt]
			\end{tikzcd}
			\caption{Curve $X$ of strict Belyi type}\label{Fig:StcB}
		\end{wrapfigure}
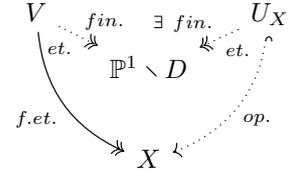
		\subsection{Combinatorial Belyi cuspidalization}
			From the perspective of algorithmic anabelian reconstruction (AlgAG), a fundamental step consists in reconstructing the cuspidal decomposition/inertia groups $I_x\simeq \widehat{\ZZ}$ from $\piet(X\times\bar{K})$. In the case where $X$ is a \emph{non-proper} hyperbolic curve, this reconstruction follows, relatively to the smooth compactification of $X$, from Nakamura's ``cyclotomic normalized'' property, as previously discussed in \S~\ref{subsub:AnabNak}.
			
			\medskip
			
			In the case of a proper hyperbolic curve, one assumes that $X$ is of \emph{strict Belyi type}, that is, that it fits in a diagram such as in Fig.~\ref{Fig:StcB}, see \cite{AbsTopII} \S~3. This condition amounts functorially relating finite étale covers to some Belyi ramification data -- see Rem.~1.11.3 \cite{AbsTopIII}. As a consequence, one can reconstruct from $\piet(X)$, as in Cor.~3.7 of \cite{AbsTopII}, the set
			
			\medskip
			
			{\centering
				$\{\piet(U_X)\longrightarrow\piet(X)\}_{U_X\leq X}$,
			\par}
			
			\medskip
			
			where $U\leq X$ runs over the open complement of finite points -- so that one can apply Nakamura's cyclotomic normalized argument of \S~\ref{subsub:AnabNak} to the various $\piet(U_X)$.
			
			\medskip
			
			The development of the anabelian properties of Galois-Teichmüller theory relies in an essential manner on the similar \emph{theory of combinatorial Belyi cuspidalization}, as introduced in \cite{Tsuj20}.

			\subsubsection{Belyi diagrams}
				A Belyi diagram is a group-theoretic setup that emanates from an open immersion/finite étale covering setup of $X_{0,3}$ as in Fig.~\ref{Diag:BelyiDiag} -- see \cite{Tsuj20} Def.~1.1~(i): \textbf{(a)}~the cusps information of the open immersion is controlled in terms of étale covers, and \textbf{(b)}~anabelian data in $\MM_{0,[m]}$ are inherited, via tripod synchronization techniques as in \cite{CbTpII} and \S~\ref{sub:BelInj}, from $\MM_{0,4}\simeq X_{0,3}$.
				
				\medskip
			
				{\centering
					$
					\begin{tikzcd}[column sep=-.7em]
						& U \ar[dr,"\text{fin. étale}",two heads] \ar[dl,"\text{op. im.}", hook', swap]& \\
						X_{0,3}=\mathbb{P}_{\bar{\QQ}}^1\setminus\{0,1,\infty\}	& &	X_{0,3}=\mathbb{P}_{\bar{\QQ}}^1\setminus\{0,1,\infty\}
					\end{tikzcd}
					\quad\rightsquigarrow\quad%
					\begin{tikzcd}[column sep=.7em]
						& \Pi_U \ar[dl,"\text{Out surj.}"',two heads] \ar[dr, "open",hook]& \\
						\Pi_{X_{0,3}}	& &	\Pi_{X_{0,3}}\overset{Out}{\reflectbox{$\lcurvearrowright$}}\, \GT
					\end{tikzcd}
					$
					
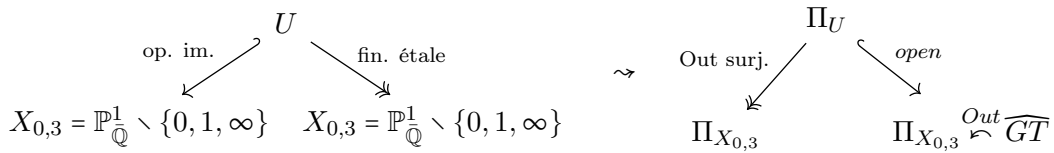
\captionof{figure}{A Belyi diagram}\label{Diag:BelyiDiag}
				\par}

			\subsubsection{Combinatorial Belyi cuspidalization} 
\begin{wrapfigure}[5]{r}{7cm}
					\vspace*{-2em}
					\centering
					$\begin{tikzcd}[column sep=1em, row sep=1em]
							& \Pi_U\rO N \ar[dl,two heads] \ar[dr,hook] \ar[d,phantom,"\mathbb{B}"]& \\
							\Pi_{X_{0,3}}\rO N	&\hphantom{A} &	\Pi_{X_{0,3}}\rO N
					\end{tikzcd}$
					\captionof{figure}{An arithmetic Belyi diagram}\label{fig:ArBelDiag}
				\end{wrapfigure}
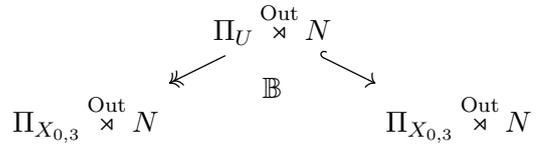

				Let us see how the $\GT$-information can be transfered from $\GT\leq \Out (\Pi_{X_{0,3}})$ to any  closed subgroup $G\leq \Out(\Pi_{X_{0,3}})$ in a combinatorial group-theoretical way. We proceed at a finite level, namely with respect to a chosen normal open subgroup $N\leq \GT$. For a closed subgroup $G\leq \Out(\Pi_{X_{0,3}})$, we denote $\Pi_{X_{0,3}}\rtimes^{\Out} G\coloneq\Aut(\Pi_{X_{0,3}})\times_{\Out(\Pi_{X_{0,3}})} G$ as in \cite{Tsuj20} -- recall that $\Pi_{X_{0,3}}$ is a center-free group. In this setup, Belyi diagram becomes \emph{arithmetic Belyi diagram} as in Fig.~\ref{fig:ArBelDiag} and ibid. Def.~1.5.
			
				\medskip
				
				The following is \cite{Tsuj20} Th.~1.3 (i) and~(iii).
				\begin{prop}[Combinatorial Belyi Cuspidalization]
					Consider given a Belyi diagram as above. Then the morphism $\Pi_U\twoheadrightarrow \Pi_{X_{0,3}}$ is group-theoretically reconstructed in a unique way from the following data
					\[
					\Pi_{X_{0,3}},\ \GT\leq \Out(\Pi_{X_{0,3}}),\ \Pi_U\leq \Pi_{X_{0,3}},\text{ and } \{0,1,\infty\}\text{ as cusps of } \Pi_{X_{0,3}}\text{ and of }\Pi_U.
					\]
          Similar results also hold for any \emph{tripodal pair} $(\Pi,G\leq \Out(\Pi))$ in the sense of ibid., that is, any pair abstractly isomorphic to $(\Pi_{X_{0,3}},\GT\leq\Out(\Pi_{X_{0,3}}))$.
				\end{prop}
				
				By finiteness of the cusps of $U$ and an Aut/Out diagram chasing between $\Pi_U$ and $\Pi_{X_{0,3}}$, one first 
				shows that any sufficiently small open normal subgroup $N\leq \GT$ identifies with a $\widetilde{N}_U\leq \Out(\Pi_U)$ in a way that

					\medskip
					
					{\centering
					$\Pi_U\rO \widetilde{N}_U\underset{open}{\hookrightarrow}\Pi_{X_{0,3}}\rO \GT,\text{ and }\widetilde{N}_U\text{ acts trivially on the cusps of } \Pi_U,$
					\par}
					
					\medskip				
	
				and so that the construction is moreover compatible with the Belyi formation (resp. the $\GT$-open formation) $\Pi_U\hookrightarrow\Pi_{X_{0,3}}$ (resp. $N\triangleleft \GT$) -- see Th.~1.3~(ii) ibid.
				
			\medskip
			
			\begin{wrapfigure}[6]{r}{7cm}
				\vspace*{-2.5em}
				\centering
				$\begin{tikzcd}[column sep=1.2em]
					&\Out^{\rm FC}(\Pi_{X_3})\ar[d, "\Pi_{X_3}\hookleftarrow\Pi_{U_3}",sloped] \ar[rd, "\Pi_{X_3}\twoheadrightarrow\Pi_{U_3}",sloped] \ar[ld,"Th.~\ref{theo:OutFTinj}"',hook']& \\[2em]
					\Out(\Pi_{X_{0,3}}) \ar[r, "\sim",dotted]& \Out(\Pi_T) \ar[r, "\sim",dotted] &  \Out(\Pi_T)
				\end{tikzcd}$
				\captionof{figure}{Tripods synchronization for $N=\widetilde{N}$ }\label{fig:TripSync3}
			\end{wrapfigure}
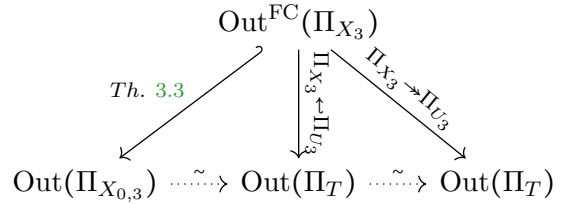
			Finally, since the cusps generate $\ker\{\Pi_U\to \Pi_{X_{0,3}}\}$, one further obtains, via $\Pi_U\twoheadrightarrow \Pi_{X_{0,3}}$, a certain subgroup $\widetilde{N}\leq \Out(\Pi_{X_{0,3}})$, which, one shows using the tripod synchronization in the section below, coincides with the original subgroup $N\triangleleft \GT\leq \Out (\Pi_{X_{0,3}})$.
			
			\subsubsection{Tripod Synchronization}

				The identification $N\simeq \widetilde{N}$ follows from the synchronization of various tripods $\Pi_T\leq \Pi_{X_3}$ (as in Fig.~\ref{fig:TripSync3}), that appears in the $3$-dimensional $\MM_{0,6}\simeq X_3$.
				The identity $ N=\widetilde{N}$ can here be read, via the injectivity property of Th.~\ref{theo:OutFTinj},  within $N\leq \Out^{\rm FC}(\Pi_{X_3})$ -- note that the right hand-side triangle of Fig.~\ref{fig:TripSync3} is given by a Belyi cuspidalization diagram.
			
				\medskip
				
				For the reader convenience, we produce, following Tsujimura, the geometric situation in Fig.~\ref{Fig:Tripods} involved in terms of tripod exceptional divisors as in \S~\ref{sub:BelInj}.

			\begin{center}
					\def\svgwidth{.5\linewidth}
					\resizebox{.45\linewidth}{!}{
\begingroup%
  \makeatletter%
  \providecommand\color[2][]{%
    \errmessage{(Inkscape) Color is used for the text in Inkscape, but the package 'color.sty' is not loaded}%
    \renewcommand\color[2][]{}%
  }%
  \providecommand\transparent[1]{%
    \errmessage{(Inkscape) Transparency is used (non-zero) for the text in Inkscape, but the package 'transparent.sty' is not loaded}%
    \renewcommand\transparent[1]{}%
  }%
  \providecommand\rotatebox[2]{#2}%
  \newcommand*\fsize{\dimexpr\f@size pt\relax}%
  \newcommand*\lineheight[1]{\fontsize{\fsize}{#1\fsize}\selectfont}%
  \ifx\svgwidth\undefined%
    \setlength{\unitlength}{85.27405686bp}%
    \ifx\svgscale\undefined%
      \relax%
    \else%
      \setlength{\unitlength}{\unitlength * \real{\svgscale}}%
    \fi%
  \else%
    \setlength{\unitlength}{\svgwidth}%
  \fi%
  \global\let\svgwidth\undefined%
  \global\let\svgscale\undefined%
  \makeatother%
  \begin{picture}(1,0.74464513)%
    \lineheight{1}%
    \setlength\tabcolsep{0pt}%
    \put(0,0){\includegraphics[width=\unitlength,page=1]{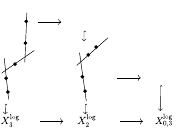}}%
    \put(0.4003068,0.61434883){\color[rgb]{0,0,0}\makebox(0,0)[lt]{\Large\smash{\begin{tabular}[t]{l}$\{\Delta\}$\end{tabular}}}}%
    \put(0.17533306,0.65669133){\color[rgb]{0,0,0}\makebox(0,0)[lt]{\smash{\begin{tabular}[t]{l}$\Delta$\end{tabular}}}}%
    \put(0.57932035,0.48008648){\color[rgb]{0,0,0}\makebox(0,0)[lt]{\smash{\begin{tabular}[t]{l}$\Delta$\end{tabular}}}}%
    \put(0.14723484,0.52527795){\color[rgb]{0,0,0}\makebox(0,0)[lt]{\smash{\begin{tabular}[t]{l}$\Delta'$\end{tabular}}}}%
    \put(0.8489247,0.30540251){\color[rgb]{0,0,0}\makebox(0,0)[lt]{\Large\smash{\begin{tabular}[t]{l}$\{\infty\}$\end{tabular}}}}%
    \put(0.06204799,0.22448783){\color[rgb]{0,0,0}\makebox(0,0)[lt]{\smash{\begin{tabular}[t]{l}$0$\end{tabular}}}}%
    \put(0.50895891,0.2250265){\color[rgb]{0,0,0}\makebox(0,0)[lt]{\smash{\begin{tabular}[t]{l}$0$\end{tabular}}}}%
    \put(0.06040521,0.304183){\color[rgb]{0,0,0}\makebox(0,0)[lt]{\smash{\begin{tabular}[t]{l}$1$\end{tabular}}}}%
    \put(0.1141979,0.38380894){\color[rgb]{0,0,0}\makebox(0,0)[lt]{\smash{\begin{tabular}[t]{l}$\infty$\end{tabular}}}}%
    \put(0.49526227,0.39089197){\color[rgb]{0,0,0}\makebox(0,0)[lt]{\smash{\begin{tabular}[t]{l}$\infty$\end{tabular}}}}%
    \put(0.48257349,0.304183){\color[rgb]{0,0,0}\makebox(0,0)[lt]{\smash{\begin{tabular}[t]{l}$1$\end{tabular}}}}%
  \end{picture}%
\endgroup%
}
					
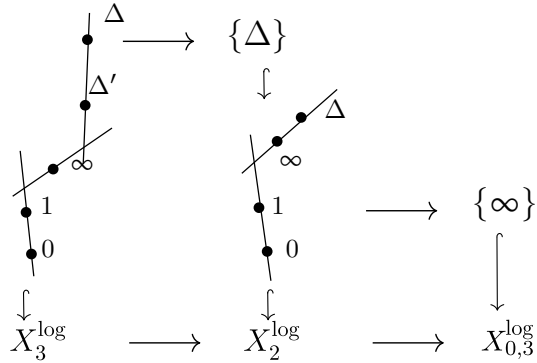
\captionof{figure}{The tripod configuration in $X_3$}\label{Fig:Tripods}
			\end{center}
				
				\medskip
				
				\begin{rmk}
					Belyi diagram can be defined in a purely combinatorial and group-theoretic manner that, in particular, does not rely on any version of the Grothendieck anabelian conjecture for hyperbolic curves -- see \cite{HMT25} Rem.~3.3.3 (and Rem.~3.3.2).
				\end{rmk}	
				
			\subsection{The BGT combinatorial models}
				Let us fix a closed subgroup $G\leq \GT$, and let us proceed at the level of some open normal subgroups $N$ and $N^\dagger\trianglelefteq G$. The obtention of combinatorial models for $G_{\QQ}$ and $\bar{\QQ}$ relies on the use of arithmetic Belyi diagrams as in Fig.~\ref{fig:ArBelDiag} to dominate $\Pi_{X_{0,3}}$ in terms of open normal subgroups $N$ of $G$.
				
				\medskip
				
				A diagram $\mathbb{B}_{N^\dagger}$ dominates $\mathbb{B}_{N}$ in the following situation
				\[\begin{tikzcd}[column sep= 1em,/tikz/column 4/.style={column sep=-.7em}]
					& \Pi_U\rO N \ar[ld,hook]\ar[rd,two heads]\ar[rrrr,two heads,dotted, "\exists \varphi\text{ on }N\cap N^\dagger","open"']\ar[d,phantom,"\mathbb{B}_N"]& & & &\Pi_{U^\dagger}\rO N^\dagger \ar[ld,hook]\ar[rd,two heads] \ar[d,phantom,"\mathbb{B}_{N^\dagger}"]& \\
					\Pi_{X_{0,3}}\rO N	& \hphantom{A} &\Pi_{X_{0,3}}\rO N &  &\Pi_{X_{0,3}}\rO N^\dagger & \hphantom{A}& \Pi_{X_{0,3}}\rO N^\dagger 
				\end{tikzcd}
				\]
				\noindent so that the $\varphi$ lifts to $\Pi_{U_2}$ to provide an outer-surjection
				\[\Pi_{U_2}\rO N\cap N^\dagger\twoheadrightarrow \Pi_{U}\rO N\cap N^\dagger\] 
				whose kernel preserves fiber subgroups and their cusps -- see \cite{HMT25} Def.~3.3~(i).
				
				\medskip
				
				\subsubsection{Belyi-Galois-Teichmüller groups} The cofiltered property (COF) and the existence (relative Grothendieck conjecture, RGC) of domination systems are, so far, encoded in the definition of a Belyi-Galois-Teichmüller group -- see \cite{HMT25} Def.~3.3~(v).
				\begin{defi}
					A BGT group is a closed subgroup $\BGT\leqslant\GT$ that is \emph{Belyi compatible}, i.e., such that every arithmetic Belyi diagrams $\mathbb{B}$ and $\mathbb{B}^\dagger$
					\begin{itemize}
						\item[(COF)] over any normal subgroup of $\BGT$ can be dominated by an arithmetic Belyi diagram dominating $\mathbb{B}$ and $\mathbb{B}^\dagger$;
						\item[(RGC)] have at most one (geometric) domination.
					\end{itemize}
				\end{defi}
				
				Consider $I_{\BGT}$ the countable set of Belyi arithmetic diagrams associated to normal open subgroups of $\BGT$. As a direct result, one obtains, from the cusps of $\Pi_U$ of $\mathbb{B}$ and by taking the direct limit over $I_{\BGT}$,  a combinatorial model $\bar{\QQ}_{\BGT}$ of the algebraic closure of the rational numbers, that is \emph{a $\BGT$-realization of $\bar{\QQ}$}  -- see \cite{HMT25} Def.~4.1:
				\[
					\bar{\QQ}_{\BGT}\coloneq\varinjlim_{\mathbb{B}\in I_{\BGT}} \operatorname{Cusp}(\mathbb{B})\setminus\{\infty\}
				\]
				
				\begin{rmk}
					A $\BGT$ Galois model (resp. a rational $\bar{\QQ}_{\BGT}$ model) is defined up to conjugacy in $\GT$ (resp. up to isomorphism).
				\end{rmk}
				
			\subsubsection{Combinatorial models} 
				We finally reach the definition of our combinatorial model of $\bar{\QQ}$; the following is Th.~4.4 of \cite{HMT25}.
				\begin{theo}
					The $\BGT$-realizations $\bar{\QQ}_{\BGT}$ of $\bar{\QQ}$ admit a field structure $(\bar{\QQ}_{\BGT},\boxtimes_{\BGT},\boxplus_{\BGT})$, and one has a field isomorphism $\bar{\QQ}\simeq \bar{\QQ}_{\BGT}$.
				\end{theo}
				
				The field structure involves furthermore two involutions
				\[
					\medsquare^{-1}\text{ and }(1-\medsquare)\colon\bar{\QQ}_{\BGT}\cup\{\infty\}\to \bar{\QQ}_{\BGT}\cup\{\infty\},
				\]
				which, to introduce a $-1$-root of unity, corresponds to a degree $2$ cover of the projective line.
				 
				 \medskip
				\begin{wrapfigure}[5]{r}{5cm}
					\centering
					\vspace{-1.6em}
					$\begin{tikzcd}[column sep=.2em]
						& \Pi_U\rO N \ar[dl,"\text{Out surj.}"',two heads] \ar[dr, "open",hook]& \\
						\Pi_{X}\rO N \ar[rr, "x\rightsquigarrow y"] & &	\Pi_{X}\rO N 
					\end{tikzcd}$
					\captionof{figure}{Definition of $\boxtimes_{\BGT}$ }
				\end{wrapfigure}
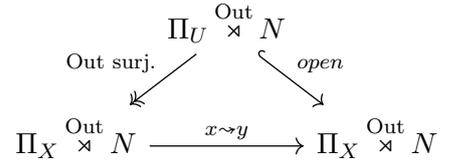
				The values of $\boxtimes_{\BGT}$ and $\boxplus_{\BGT}$ are defined using function field-like arguments -- in the sense that $\mathbb{P}^1_{\bar{\QQ}}\setminus\{0,1,\infty\}\simeq {\rm Spec}\,\bar{\QQ}[t,1-t,1/t]$ -- from the cusp properties of Belyi diagrams as follows.
				
				In a Belyi diagram $\mathbb{B}$, replacing $\{0,1,\infty\}\subset\operatorname{Cusp}(\mathbb{B})= \operatorname{Cusp}(U\rtimes^{\Out} N)$ by $\{\infty,1,0\}$, $\{1,0,\infty\}$, and $\{0,x,\infty\}$ for $x\in \bar{\QQ}_{\BGT}$ defines different values $y$, and, by cofinality, some functions:
				\[
					\{t^{-1}\}\colon \bar{\QQ}_{\BGT}\setminus\{0\} \to \bar{\QQ}_{\BGT}\setminus\{0\}, \text{ and } \{1-t\},\{t/x\}\colon \bar{\QQ}_{\BGT}\to \bar{\QQ}_{\BGT},
				\]
				which, for $x,\, y\in \bar{\QQ}_{\BGT}$, define the different values:
				\[
				\boxtimes_{\BGT}(x,y)\coloneq\{t/\{t^{-1}(x)\}(y)\}\quad\text{ with }
				\begin{cases}
					\boxtimes_{\BGT}(0,y)\coloneq 0\\
					\boxtimes_{\BGT}(1,y)\coloneq y
				\end{cases}
				\]
				 The values of $\boxplus_{\BGT}$ are defined in a similar but more intricate manner -- see ibid. before Claim 4.4~A and Claim 4.4~D.
				
				\medskip
				
				\begin{wrapfigure}[10]{r}{5cm}
					\centering
					\vspace{-.9em}
					$\begin{tikzcd}[column sep=.2em]
						&\GT & & \\[-.5em]
						G_\QQ^\sigma\ar[ru,phantom,"\leq",sloped] \ar[d,rightsquigarrow] & & G_\QQ^\tau \ar[lu,phantom,"\leq",sloped] \ar[d,rightsquigarrow] & \BGT\\
						\prescript{\dagger}{}{\bar{\QQ}}^\sigma_{\BGT} \ar[rr,dotted,"\sigma.\tau^{-1}"']& & \prescript{\ddagger}{}{\bar{\QQ}}^\tau_{\BGT} & \bar{\QQ}_{\BGT}
					\end{tikzcd}$
					\captionof{figure}{$\BGT$ rational models of $\bar{\QQ}$}\label{fig:BGTmod}
				\end{wrapfigure}
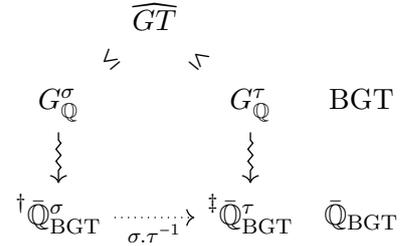
				The isomorphism with the algebraic closure of the rational numbers follows from a similar group-theoretic construction of any finite subset of $\bar{\QQ}$ in terms of Belyi diagrams to obtain a field homomorphism $\bar{\QQ}\hookrightarrow \bar{\QQ}_{\BGT}$ which is an isomorphism by construction -- see \cite{HMT25} Th.~4.4~(ii).

				\subsubsection{BGT and anabelian variation of holomorphic structure} In the global construction, up-to-conjugacy, of some $\BGT$ models of $\bar{\QQ}$ from (a closed subgroup $G$ of) $\GT$, it is of importance to note that the mono-anabelian algorithms into play involve \emph{multiple distinct ring structures} -- see Fig.~\ref{fig:BGTmod}. As a consequence, \textbf{(a)}~the isomorphism $\sigma.\tau^{-1}$ between $\BGT$-models is not algebraic, and it can be seen as the result of a \emph{variation of analytic structure}\footnote{Since the natural outputs of the anabelian algorithms involved are multiplicative monoids, such a variation of analytic structure can be understood, at the level of ring structures, as a variation of the additive monoid structure.} -- as for the parameter $f_\sigma\in\widehat{\mathbb{F}}_2'$ of \S~\ref{sub:P3pts} --  and \textbf{(b)}~ the group $\GT$ (or the closed subgroup $G$) can be seen as \emph{the automorphism group of an arithmetic-holomorphic structure} -- see (iii) in \S~\ref{subsub:AGprinc}.

					\begin{rmk}\label{rm:freeGT}
						The constructions of this section can be realized without any reference to the Grothendieck-Teichmüller group\footnote{Thus the interpretation in this paper, which was left to the reader in the original \cite{HMT25}, of ``BGT'' as a ``Belyi Galois-Teichmüller'' model.}, but to the anabelian properties of Galois-Teichmüller theory only, that is of $\Out(\Pi_{X_{0,3}})$.
					\end{rmk}
				
				\subsection{Opening perspectives...}
					Given the original multifaceted geometric and arithmetic nature of Grothendieck-Teichmüller theory  -- in relation with braid groups, mapping class groups of diffeomorphisms, or Galois properties -- and the recent universal anabelian approach described in this manuscript (see also (RNS) below), we gather some topics of potential future interest.
					
					\subsubsection{Arithmetic, Braids, mapping class groups}
						
						In the original \emph{braid group approach} of \S~\ref{sub:higher}, the Galois-Teichmüller anabelian isomorphism of Eq.~\eqref{eq:GTanab} identifies the outer automorphism groups of certain braid group quotient and mapping class groups, see \cite{MN22}:
						\[
							\Out(\widehat{B}_m/Z_m)\simeq\GT \text{ and } \Out(\widehat{\Gamma}_{1,2})\simeq \GT
						\]
						where $Z_m$ denotes the center of $\widehat{B}_m$.
						
						\medskip
						
						The \textit{strong indecomposability} property -- that is for a profinite group, up to open subgroup, to admit only trivial direct decomposition -- was first established \textbf{(a)}~for absolute Galois group of number fields \cite{HJ91} Cor.~2.5, to turn into an anabelian property satisfied by \textbf{(b)}~the étale fundamental group of surface groups, of hyperbolic curves over number fields (and configuration spaces of, and outer automorphisms of), and \textbf{(c)}~also $\GT$ -- see \cite{MTam08} and \cite{MT25}. 
						
						It now has motivated a generalization of Jarden's ``families preserving'' isomorphisms property for Galois groups in relation with their inner automorphisms, see \cite{MST24}.
					
					\subsubsection{Combinatorial homotopy over $p$-adic local fields}
						Over $K$ a $p$-adic local field, the $p$-adic topology and successive combinatorial considerations each produces various $p$-adic avatars of the Grothendieck-Teichmüller group $\GT$
						\[
							G_{\QQ_p}\leq {\rm GT}^{\rm M},\ {\rm GT}_p\leq {\rm GT}_p^{\rm tp}\coloneq \GT\cap \Out(\Pi^{\rm tp}_{X_{0,3}})\leq \Out(\Pi_{X_{0,3}})
						\]
						where ${\rm GT}_p$ denotes André's $p$-adic $\GT$ avatar of \cite{And03} that results of tower considerations as in \S~\ref{sub:GTtower}, ${\rm GT}^{\rm M}$ is the metrized version of $\GT$ of \cite{CbTpIII}, ${\rm GT}_p^{\rm tp}$ is obtained in \cite{Tsuj20}, and the homotopy context is given by André's tempered fundamental group $\pi_1^{\rm tp}(X)$ -- see  \cite{Lep10} for a functorial approach and references herein.

						\medskip
						
						The following is a consequence of establishing Tamagawa-Lepage's ``resolution of nonsingularity'' (RNS) for hyperbolic curve over $p$-adic local fields, see \cite{RNS23} Th.~G:
						\[
						 \text{\itshape The following equalities hold }{\rm GT}^{\rm M}= {\rm GT}_p={\rm GT}_p^{\rm tp}\leq \Out(\Pi_{X_{0,3}})
						\]
						From the point of view of anabelian geometry, a direct consequence of (RNS) is to establish the \emph{absolute} Grothendieck's conjecture (GC) for hyperbolic curves (resp. configuration spaces of) over $p$-adic local fields, see ibid. Th.~D (resp. Th.~E).
						
						\medskip
						
						The following is a direct consequence of the BGT approach, see \cite{HMT25} Rem.~4.4.1
						\begin{prop}\label{prop:GTGF} 
							For $K$ a $p$-adic local field, after identification of $\bar{\QQ}$ with the algebraic closure of $\QQ$ in $\bar{K}$, the Galois group $G_K$ identifies as a subgroup of $\GT$ and satisfies both the (COF) and the (RGC) properties, so that $G_K\simeq$ ``$\BGT$'' and $\bar{\QQ}_{G_K}\simeq\bar{\QQ}$.
						\end{prop}
						
					\bigskip

						\subsubsection{Towards a definite combinatorial model of $G_\QQ$}
							In light of the universality inherent in anabelian arithmetic geometry within Galois-Teichmüller theory -- in the senses that homotopy essentially and faithfully encodes both the arithmetic and the geometry of the spaces involved -- together with the decisive progress of \S~\ref{sec:GTAnab} and \S~\ref{sec:BGT}, it now appears natural to investigate the definition of a \emph{new, purely anabelian combinatorial model ${\rm GT}_{\rm NF}$}, as hinted in Rem.~\ref{rm:freeGT}, fitting into:			 
							\[
								G_\QQ \leq {\rm GT}_{\rm NF} \leq \widehat{GT}\leq \Out(\Pi_m).
							\]
							
							\medskip
							
							Defined using suitable classes of Belyi maps so as to reflect the intrinsic ``genus zero'' nature of $\GT$, and in analogy with the (COF) $p$-adic local result of Prop.~\ref{prop:GTGF}, this group would constitute a compelling candidate for an isomorphism with the absolute Galois group $G_\QQ$. We refer to the forthcoming \cite{ArGT26}.



\newpage
			
{\scshape\raggedleft\large References\par}
\addcontentsline{toc}{section}{References}
\begin{multicols}{2}
	\printbibliography[heading=none] 
\end{multicols}

\end{document}